\input ./style/arxiv-general.cfg
\documentclass[MSNbibl,number,citesort,seceqn,dvips]{arxbj}
\makeatletter
   \@ifpackageloaded{graphicx}{}{\usepackage{graphicx}}
\makeatother

\volume{22}
\issue{4}
\pubyear{2016}
\firstpage{2001}
\lastpage{2028}
\doi{10.3150/15-BEJ718}
\docsubty{FLA}

\makeatletter
\newcommand{\rrvert}{\vert}
\newcommand{\llvert}{\vert}
\newremark{problem}{Problem}
\newproclaim{cond}{Condition}
\newtheorem{theorem}{Theorem}[section]
\newtheorem{lemma}[theorem]{Lemma}
\newtheorem{proposition}[theorem]{Proposition}
\newproclaim{definition}[theorem]{Definition}
\newtheorem{corollary}[theorem]{Corollary}
\newproclaim{assumption}[theorem]{Assumption}
\newremark{example}[theorem]{Example}
\newremark{remark}[theorem]{Remark}
\makeatother

\begin{document}
\begin{frontmatter}

\title{Quantifying repulsiveness of determinantal point processes}
\runtitle{Quantifying repulsiveness of determinantal point processes}

\begin{aug}
\author[A]{\inits{C.A.N.}\fnms{Christophe Ange Napol\'eon}~\snm
{Biscio}\corref{}\thanksref{A}\ead[label=e1]{Christophe.Biscio@univ-nantes.fr}}
\and
\author[B]{\inits{F.}\fnms{Fr\'ed\'eric}~\snm{Lavancier}\thanksref
{B}\ead[label=e2]{Frederic.Lavancier@univ-nantes.fr}}
\address[A]{Laboratoire de Math\'ematiques Jean Leray -- BP 92208 --
2, Rue de la Houssini\`ere -- F-44322 Nantes Cedex 03 -- France.
\printead{e1}}
\address[B]{Inria, Centre Rennes Bretagne Atlantique, France.
\printead{e2}}
\end{aug}

%
\received{\smonth{7} \syear{2014}}
%
\revised{\smonth{2} \syear{2015}}

%
\begin{abstract}
Determinantal point processes (DPPs) have recently proved to be a
useful class of models in several areas of statistics, including
spatial statistics, statistical learning and telecommunications
networks. They are models for repulsive (or regular, or inhibitive)
point processes, in the sense that nearby points of the process tend to
repel each other.
We consider two ways to quantify the repulsiveness of a point process,
both based on its second-order properties, and we address the question
of how repulsive a stationary DPP can be.
We determine the most repulsive stationary DPP, when the intensity is
fixed, and for a given $R>0$ we investigate repulsiveness in the
subclass of $R$-dependent stationary DPPs, that is, stationary DPPs
with $R$-compactly supported kernels. Finally, in both the general case
and the $R$-dependent case, we present some new parametric families of
stationary DPPs that can cover a large range of DPPs, from the
stationary Poisson process (the case of no interaction) to the most
repulsive DPP.
\end{abstract}

%
\end{frontmatter}

\section{Introduction}\label{sectionintroduction}

Determinantal point processes (DPPs) were introduced in their general
form by Macchi \cite{macchi1975coincidence} in 1975 to model fermions in
quantum mechanics, though some specific DPPs appeared much earlier in
random matrix theory. DPPs actually arise in many fields of probability
and have deserved a lot of attention from a theoretical point of view,
see for instance \cite{hough2009zeros} and \cite{Soshnikov00}.

DPPs are repulsive (or regular, or inhibitive) point processes, meaning
that nearby points of the process tend to repel each other (this
concept will be clearly described in the following). This property is
adapted to many statistical problems where DPPs have been recently
used, for instance in telecommunication to model the locations of
network nodes \cite{deng2014ginibre,miyoshi2014} and in statistical
learning to construct a dictionary of diverse sets
\cite{KuleszaTaskar12}. Other examples arising from biology, ecology and
forestry are studied in \cite{lavancierpublish} and its associated
on-line supplementary file \cite{lavancierextended}.

The growing interest for DPPs in the statistical community is due to
that their moments are explicitly known, parametric families can easily
been considered, their density on any compact set admits a closed form
expression making likelihood inference feasible and they can be
simulated easily and quickly. Section~\ref{sectionstationaryDPP}
summarizes some of these properties and we refer to \cite
{lavancierpublish} for a detailed presentation.
These features make the class of DPPs a competitive alternative to the
usual class of models for repulsiveness, namely the Gibbs point
processes. In contrast, for Gibbs point processes, no closed form
expression is available for the moments, the likelihood involves an
intractable normalizing constant and their simulation requires Markov
chain Monte Carlo methods.

However, DPPs cannot model all kinds of repulsive point patterns. For
instance, as deduced from Section~\ref{sectionDPPregularitypcf},
stationary DPPs cannot involve a hardcore distance between points,
contrary to the Mat\'ern hardcore point processes, the RSA (random
sequential absorption) model and hardcore Gibbs models; see \cite
{illianpenttinenstoyanstoyan08}, Section~6.5. In this paper, we
address the question of how repulsive a stationary DPP can be. We also
investigate for a given $R>0$ the repulsiveness in the subclass of
$R$-dependent stationary DPPs, that is, stationary DPPs with
$R$-compactly supported kernels, which are of special interest for
statistical inference in high dimension, see Section~\ref{sectionDPPregularitypcffiniterange}. In both cases, we present in
Section~\ref{DPPkernelsfamily} some parametric families of stationary DPPs that
cover a~large range of DPPs, from the stationary Poisson process to the
most repulsive DPP.

To quantify the repulsiveness of a stationary point process, we
consider its second-order properties. Let $X$ be a stationary point
process in $\mathbb R^d$ with intensity\vspace*{1pt} (i.e. expected number of points
per unit volume) $\rho>0$ and second order intensity function $\rho
^{(2)}(x,y)$. Denoting $dx$ an infinitesimal region around $x$ and
$|dx|$ its Lebesgue measure,\vspace*{1pt} $\rho|dx|$ may be interpreted as the
probability that $X$ has a point in $dx$. For $x\neq y$, $\rho
^{(2)}(x,y)|dx||dy|$ may be viewed as the probability that $X$ has a
point in $dx$ and another point in $dy$. A formal definition is given
in Section~\ref{sectionstationaryDPP}. Note that $\rho
^{(2)}(x,y)=\rho
^{(2)}(0,y-x)$ is a symmetric function and depends only on $y-x$
because of our stationarity assumption.

In spatial statistics, the second-order properties of $X$ are generally
studied through the pair correlation function (in short p.c.f.), defined
for any $x\in\mathbb{R}^d$ by
\[
g(x)=\frac{\rho^{(2)}(0,x)}{\rho^2}.
\]
Since $\rho^{(2)}$ is unique up to a set of Lebesgue measure zero (see
\cite{daleyvol1}), so is $g$. As it is implicitly done in the
literature (see \cite
{illianpenttinenstoyanstoyan08,stoyan1987stochastic}), we choose
the version of $g$ with as few discontinuity points as possible.
It is commonly accepted (see, e.g., \cite{stoyan1987stochastic}) that
if $g(x)=1$ then there is no interaction between two points separated
by $x$, whereas there is attraction if $g(x)>1$ and repulsiveness if
$g(x)<1$. Therefore, when we below compare the global repulsiveness of
two stationary point processes, we assume they share the same intensity.

\begin{definition}\label{globallyrepulsive}
Let $X$ and $Y$ be two stationary point processes with the same
intensity $\rho$ and respective pair correlation function $g_X$ and
$g_Y$. Assuming that both $(1-g_X)$ and $(1-g_Y)$ are integrable, we
say that $X$ is globally more repulsive than $Y$ if $\int(1-g_X) \geq
\int(1-g_Y)$.
\end{definition}

The quantity $\int(1- g)$ is already considered in the on-line
supplementary material \cite{lavancierextended} of \cite
{lavancierpublish} as a measure for repulsiveness. It can be justified
in several ways. First, it is a natural geometrical method to quantify
the distance from $g$ to 1 (corresponding to no interaction), where the
area between $g$ and $1$ contributes positively to the measure of
repulsiveness when $g<1$ and negatively if $g>1$. Second, denoting $K$
and $K_0$ the Ripley's $K$-functions of $X$ and of the stationary
Poisson process with intensity $\rho$ respectively (see \cite
{mollerstatisticalinference}, Definition~4.6), we have $\int(1-g) =
\lim_{r\to\infty} (K_0 (r) - K(r))$. We also refer to \cite
{lavancierextended} for an equivalent interpretation in terms of the
reduced Palm distribution. Finally, it is worth mentioning that for any
stationary point processes, we have $\int(1- g)\leq1/\rho$, see
\cite{kuna07}, equation (2.5).

Additional criteria could be introduced to quantify the global
repulsiveness of a point process, relying for instance on $\int(1-g)^p$
for a given $p>0$, or involving higher moments of the point process
through the joint intensities of order $k>2$ (see Definition~\ref
{jointintensity}). However, the theoretical study becomes more
challenging in
these cases and we do not consider these extensions.

Repulsiveness is often interpreted in a local sense:
This is the case for hardcore point processes, where a minimal distance
$\delta$ is imposed between points and so $g(x)=0$ whenever
$|x|<\delta
$ where for a vector $x$, $|x|$ denotes its Euclidean norm. As already
mentioned, a DPP cannot involve any hardcore distance, but we may want
its p.c.f. to satisfy $g(0)=0$ and stay as close as possible to 0 near the
origin. This leads to the following criteria to compare the \textit{local
repulsiveness} of two point processes. We denote by $\nabla g$ and
$\Delta g$ the gradient and the Laplacian of $g$, respectively.

\begin{definition}\label{locallyrepulsive}
Let $X$ and $Y$ be two stationary point processes with the same
intensity $\rho$ and respective pair correlation function $g_X$ and
$g_Y$. Assuming that $g_X$ is twice differentiable at $0$ with
$g_X(0)=0$, we say that $X$ is more locally repulsive than $Y$ if
either $g_Y(0)>0$, or
$g_Y$ is not twice differentiable at $0$, or $g_Y$ is twice
differentiable at $0$ with $g_Y(0)=0$ and $\Delta g_Y(0) \geq\Delta g_X(0)$.
\end{definition}

As suggested by this definition, a stationary point process is said to
be locally repulsive if its p.c.f. is twice differentiable at $0$ with
$g(0)=0$. In this case $\nabla g(0)=0$ because $g(x)=g(-x)$. Therefore,
to compare the behavior of two such p.c.f.s near the origin, specifically
the curvatures of their graphs near the origin, the Laplacian operator
is involved in Definition~\ref{locallyrepulsive}. As an example, a
stationary hardcore process is locally more repulsive than any other
stationary point process because in this case $g(0)=0$ and $\Delta g(0)=0$.

\begin{figure}

\includegraphics{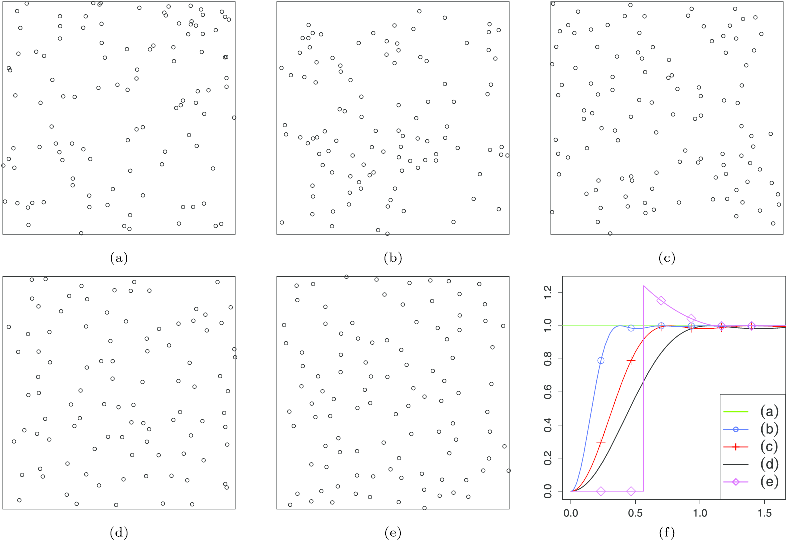}

\caption{Realizations\vspace*{1pt} on $[-5,5]^2$ of \textup{(a)} the stationary Poisson
process, \textup{(b)}--\textup{(d)} DPPs with kernels (\protect\ref
{bessel-type}) where $\sigma=0$
and $\alpha=0.2,0.4,\frac{1}{\sqrt\pi}$, \textup{(e)} the type II
Mat\'ern
hardcore process with hardcore radius $\frac{1}{\sqrt{\pi}}$.
\textup{(f)}~Their
associated theoretical p.c.f.s. The intensity is $\rho=1$ for all models
and \textup{(d)} represents the most repulsive stationary DPP in this
case.}\label{figureoverview}
\end{figure}

We show in Section~\ref{sectionDPPregularitypcf} that Definitions~\ref
{globallyrepulsive} and \ref{locallyrepulsive} agree for the natural
choice of what can be considered as the most repulsive DPP. A
realization of the latter on $[-5,5]^2$ is represented
in Figure~\ref{figureoverview}(d) when $\rho=1$. For comparison,
letting $\rho=1$
for all plots, Figure~\ref{figureoverview} shows realizations of: (a)
the stationary Poisson process, which is a situation with no
interaction; (b)--(c) two DPPs with intermediate repulsiveness, namely
DPPs with kernels (\ref{bessel-type}) where $\sigma=0$ and $\alpha
=0.2,0.4$, respectively, as presented in Section~\ref{besselfamily};
(e) the type II Mat\'ern hardcore process with hardcore radius $\frac
{1}{\sqrt{\pi}}$.
Notice that $\frac{1}{\sqrt{\pi}}$ is the maximal hardcore radius that
a type II Mat\'ern hardcore process with unit intensity can reach; see
\cite{illianpenttinenstoyanstoyan08}, Section~6.5. It corresponds
to an infinite intensity of the underlying Poisson process and our
simulation is only an approximation. These models are sorted from (a)
to (e) by their ascending repulsiveness in the sense of Definition~\ref
{locallyrepulsive}. Specifically, \mbox{$g(0) = 1$} for (a) while $g(0)=0$
and $\Delta g(0)$ is $50,  12.5,  2\pi$ and $0$ from (b) to (e),
respectively. This order is clearly apparent in~Figure~\ref{figureoverview}(f), where the theoretical p.c.f.s are represented as radial
functions, all aforementioned models being isotropic. Concerning global
repulsiveness, we have that $\int(1-g)$ is $0, 0.12, 0.50, 1$ and
$0.76$ from (a) to (e), respectively. The fact that the Mat\'ern
hardcore model is globally less repulsive than the DPP in (d) is due to
that its p.c.f. can be larger than one. This shows the limitation of
Definition~\ref{globallyrepulsive} in the study of repulsiveness and
the importance of introducing Definition~\ref{locallyrepulsive}.
Overall, Figure~\ref{figureoverview} illustrates that even if
stationary DPPs cannot be as (locally) repulsive as hardcore point
processes, which may be an important limitation in practice, they
nonetheless cover a rather large variety of repulsiveness from~(a) to~(d).

We recall the definition of a stationary DPP and some related basic
results in~Section~\ref{sectionstationaryDPP}. Section~\ref{sectionDPPregularitypcf} is devoted to the study of repulsiveness in
stationary DPPs, both in the sense of Definitions~\ref{globallyrepulsive} and~\ref{locallyrepulsive}. In
Section~\ref{sectionDPPregularitypcffiniterange}, we focus on repulsiveness
for the subclass of stationary DPPs with compactly supported kernels.
Then, in Section~\ref{DPPkernelsfamily}, we present three parametric
families of DPPs which cover a large range of repulsiveness and have
further interesting properties. Section~\ref{proofs} gathers the proofs
of our theoretical results. Further comments and illustrations are
provided in the supplementary material \cite{bisciosupplementary}.

\section{Stationary DPPs}\label{sectionstationaryDPP}

In this section, we review the definition and some properties of
stationary DPPs. For a detailed presentation, including the
nonstationary case, we refer to the survey by Hough et~al. \cite
{hough2009zeros}.

Basics of point processes may be found in \cite{daleyvol1,daleyvol2}.
Let us recall that a point process $X$ is simple if two points of $X$
never coincide, almost surely. The joint intensities of $X$ are defined
as follows.

\begin{definition}\label{jointintensity}
If it exists, the joint intensity of order $k$ ($k\geq1$) of a simple
point process $X$ is the function $\rho^{(k)}: (\mathbb{R}^d)^k \to
\mathbb{R}^+$ such
that for any family of mutually disjoint subsets $D_1,\dots,D_k$ in~$\mathbb{R}^d$,
\[
E \prod_{i=1}^k X(D_i) =
\int_{D_1}\cdots\int_{D_k} \rho
^{(k)}(x_1,\dots,x_k) \,dx_1\cdots
\,dx_k,
\]
where $X(D)$ denotes the number of points of $X$ in $D$ and $E$ is the
expectation over the distribution of $X$.
\end{definition}

In the stationary case, $ \rho^{(k)}(x_1,\ldots,x_k)= \rho
^{(k)}(0,x_2-x_1,\ldots,x_k-x_1)$, so that the intensity $\rho$ and the
second-order intensity function $\rho^{(2)}$ introduced previously
become the particular cases associated to $k=1$ and $k=2$, respectively.

\begin{definition}\label{DPPdefinition}
Let $C: \mathbb{R}^d \rightarrow\mathbb{R}$ be a function.
A point process $X$ on $\mathbb{R}^d$ is a stationary DPP with kernel
$C$, in
short $X\sim\operatorname{DPP}(C)$, if
for all $k\geq1$ its joint intensity of order $k$ satisfies the relation
\[
\rho^{(k)}(x_1,\ldots x_k)=
\det[C](x_1,\dots,x_k)
\]
for almost every $(x_1,\dots,x_k)\in(\mathbb{R}^d)^k$, where
$[C](x_1,\dots
,x_k)$ denotes the matrix with entries $C(x_i-x_j)$, $1\leq i,j\leq k$.
\end{definition}

It is actually possible to consider a complex-valued kernel $C$, but
for simplicity we restrict ourselves to the real case. A first example
of stationary DPP is the stationary Poisson process with intensity
$\rho
$. It corresponds to the kernel
\begin{equation}
\label{kernelpoisson} C(x)=\rho\mathbf{1}_{\lbrace x=0 \rbrace}\qquad \forall x\in
\mathbb{R}^d.
\end{equation}
However, this example is very particular and represents in some sense
the extreme case of a DPP without any interaction, while DPPs are in
general repulsive as discussed at the end of this section.

Definition~\ref{DPPdefinition} does not ensure existence or unicity of
$\operatorname{DPP}(C)$, but if it exists, then it is unique; see
\cite{hough2009zeros}.
Concerning existence, a general result, including the nonstationary
case, was proved by Macchi \cite{macchi1975coincidence}. It relies on the
Mercer representation of $C$ on any compact set. Unfortunately, this
representation is known only in a few cases, making the conditions
impossible to verify in practice for most functions $C$. Nevertheless,
the situation becomes simpler in our stationary framework, where the
conditions only involve the Fourier transform of $C$.
We define the Fourier transform of a function $h\in L^1(\mathbb{R}^d)$ as
\begin{equation}
\label{definitionfouriertransform} \mathcal{F}(h) (t)=\int_{\mathbb{R}^d} h(x)
e^{2i\pi x\cdot t}\,dx \qquad\forall t \in\mathbb{R}^d.
\end{equation}
By Plancherel's theorem, this definition is extended to $L^2(\mathbb{R}^d)$;
see \cite{stein1971fourier}. If $C$ is a covariance function, as
assumed in the following, we have $\mathcal{F}\mathcal{F}(C)=C$ so
$\mathcal{F}^{-1}=\mathcal{F}$ and from
\cite{sasvari2013multivariate}, Theorem~1.8.13, $\mathcal{F}(C)$
belongs to
$L^1(\mathbb{R}^d)$.

\begin{proposition}[(Lavancier, M{\o}ller and Rubak \protect\cite
{lavancierpublish})]\label{DPPexistence}
Assume $C$ is a symmetric continuous real-valued function in
$L^2(\mathbb{R}
^d)$. Then $\operatorname{DPP}(C)$ exists if and only if $0\leq
\mathcal{F}(C)\leq1$.
\end{proposition}

In other words, Proposition~\ref{DPPexistence} ensures existence of
$\operatorname{DPP}(C)$ if $C$ is a continuous real-valued covariance
function in
$L^2(\mathbb{R}^d)$ with $\mathcal{F}(C)\leq1$. Henceforth, we
assume the following
condition.

\renewcommand{\thecond}{$\mathcal{K}(\rho)$}
\begin{cond}
A kernel $C$ is said to verify condition
$\mathcal{K}(\rho)$ if $C$ is a symmetric continuous real-valued
function in $ L^2(\mathbb{R}
^d)$ with $C(0)=\rho$ and $0\leq\mathcal{F}(C)\leq1$.
\end{cond}

The assumption $0\leq\mathcal{F}(C)\leq1$ is in accordance with
Proposition~\ref{DPPexistence}, while the others assumptions in
condition $\mathcal{K}(\rho)$ are satisfied by most statistical
models of covariance
functions, the main counterexample being~(\ref{kernelpoisson}).
Standard parametric families of kernels include the Gaussian, the
Whittle--Mat\'ern and the generalized Cauchy covariance functions,
where the condition $\mathcal{F}(C)\leq1$ implies some restriction on the
parameter space; see \cite{lavancierpublish}.

By Definition~\ref{DPPdefinition}, all moments of a DPP are explicitly
known. In particular, assuming condition $\mathcal{K}(\rho)$, the
intensity of
$\operatorname{DPP}(C)$ is $\rho$ and denoting $g$ its p.c.f. we have
\begin{equation}
\label{DPPpcf} 1-g(x)=\frac{C(x)^2}{\rho^2}
\end{equation}
for almost every $x\in\mathbb{R}^d$. Consequently, $g\leq1$, and so
we have
repulsiveness.
Moreover, the study of repulsiveness of stationary DPPs, as defined in
Definitions~\ref{globallyrepulsive} and~\ref{locallyrepulsive}, reduces
to considerations on the kernel $C$ when condition $\mathcal{K}(\rho
)$ is assumed.

\section{Most repulsive DPPs}\label{sectionDPPregularitypcf}

We first present the most globally repulsive DPPs in the sense of
Definition~\ref{globallyrepulsive}. They are introduced in the on-line
supplementary file associated to \cite{lavancierpublish} (see \cite
{lavancierextended}), from which the following proposition is easily deduced.

\begin{proposition}[(Lavancier, M{\o}ller and Rubak \cite
{lavancierpublish})]\label{jinc}
In the sense of Definition~\ref{globallyrepulsive}, $\operatorname
{DPP}(C)$ is the
most globally repulsive DPP among all DPPs with kernel satisfying
condition $\mathcal{K}(\rho)$ if and only if $\mathcal{F}(C)$ is
even and equals almost
everywhere an indicator function of a Borel set with volume $\rho$.
\end{proposition}

According to Proposition~\ref{jinc}, the set of the most globally
repulsive DPPs in the sense of Definition~\ref{globallyrepulsive} is
infinite. This is illustrated in the supplementary material~\cite
{bisciosupplementary}.
A natural choice is $\operatorname{DPP}(C_B)$ where $\mathcal
{F}(C_B)$ is the indicator
function of the Euclidean ball centered at $0$ with volume~$\rho$. In
dimension $d$, this gives $C_{B}=\mathcal{F} (\mathbf{1}_{
\lbrace|\cdot|^d
\leq\rho\tau^d  \rbrace} )$ with $\tau= \lbrace
\Gamma
(d/2+1)/\pi^{{d}/{2}} \rbrace^{{1}/{d}}$ and by \cite
{grafakos2008classical}, Appendix~B.5,
\begin{equation}
\label{expressionjinc} C_{B}(x)= \frac{\sqrt{\rho\Gamma({d}/{2}+1)}}{\pi^{d/4}} \frac
{J_{{d}/{2}} (2\sqrt{\pi}\Gamma({d}/{2}+1)^{{1}/{d}}\rho
^{{1}/{d}} |x| )}{|x|^{{d}/{2}}}\qquad
\forall x\in\mathbb{R}^d,
\end{equation}
where $J_{{d}/{2}}$ is the Bessel function of the first kind. For
example, we
have
\begin{itemize}
\item for $d=1$, $C_{B}(x)=\operatorname{sinc}(x)=\frac{\sin(\pi
\rho
|x|)}{\pi|x|}$,
\item for $d=2$, $C_{B}(x)=\operatorname{jinc}(x)=\sqrt{\rho}\frac
{J_1(2\sqrt
{\pi\rho}|x|)}{\sqrt{\pi}|x|}$.
\end{itemize}

This choice was already favored in \cite{lavancierpublish}. However,
there is no indication from Proposition~\ref{jinc} to suggest $C_{B}$
instead of another kernel given by the proposition.
This choice becomes clear if we look at the local repulsiveness as
defined in Definition~\ref{locallyrepulsive}.

\begin{proposition}\label{maxlocrep}
In the sense of Definition~\ref{locallyrepulsive}, the most locally repulsive
DPP among all DPPs with kernel satisfying condition $\mathcal{K}(\rho
)$ is
$\operatorname{DPP}(C_{B})$.
\end{proposition}

Thus, from Propositions~\ref{jinc} and~\ref{maxlocrep}, we deduce the
following corollary.

\begin{corollary}\label{mostrepulsiveDPP}
The kernel $C_B$ is the unique kernel $C$ verifying condition $\mathcal
{K}(\rho)$
such that $\operatorname{DPP}(C)$ is both the most
globally and the most locally repulsive DPP among all stationary DPPs
with intensity $\rho>0$.
\end{corollary}

Borodin and Serfaty in \cite{serfaty} characterize in dimension $d\leq
2$ the disorder of a point process by its ``renormalized energy''. In
fact, the smaller the renormalized energy, the more repulsive the point
process. Theorem~3 in~\cite{serfaty} establishes that $\operatorname
{DPP}(C_B)$
minimizes the renormalized energy among all stationary DPPs. This
result confirms Corollary~\ref{mostrepulsiveDPP}, that the most
repulsive stationary DPP, if any has to be chosen, is $\operatorname
{DPP}(C_B)$.
However, a stationary DPP has a finite renormalized energy if and only
if it is given by Proposition~\ref{jinc} (\cite{serfaty}, Theorem~1),
which indicates that most stationary DPPs have an infinite renormalized
energy. Hence, this criteria is not of practical use to compare the
repulsiveness between two arbitrary DPPs.

\section{Most repulsive DPPs with compactly supported kernels}\label
{sectionDPPregularitypcffiniterange}

In this section, we assume that the kernel $C$ is compactly supported,
that is, there exists $R>0$ such that $C(x)=0$ if $|x|>R$.
In this case, $X\sim\operatorname{DPP}(C)$ is an $R$-dependent point
process in the
sense that if $A$ and $B$ are two Borel sets in $\mathbb{R}^d$
separated by a
distance larger than $R$, then $X\cap A$ and $X\cap B$ are independent,
which is easily verified using Definition~\ref{DPPdefinition}.
This situation can be particularly interesting for likelihood inference
in presence of a large number of points. Assume we observe $\{x_1,\dots
,x_n\}$ on a compact window $W\subset\mathbb R^d$, then the likelihood
is proportional to $\det[\tilde C](x_1,\ldots,x_n)$ where $\tilde C$
expresses in terms of $C$ and inherits the compactly supported property
of $C$; see \cite{lavancierpublish,macchi1975coincidence}. While this
determinant is expensive to compute if $\tilde C$ is not compactly
supported and $n$ is large, the situation becomes more convenient in
the compactly supported case, since $[\tilde C](x_1,\ldots,x_n)$ is
sparse when $R$ is small with respect to the size of $W$.
We are thus interested in DPPs with kernels satisfying the following
condition.

\renewcommand{\thecond}{$\mathcal{K}_c(\rho,R)$}
\begin{cond}
A kernel $C$ or $\operatorname{DPP}(C)$ is said to verify condition
$\mathcal{K}_c(\rho,R)$ if
$C$ verifies condition $\mathcal{K}(\rho)$ and $C$ is compactly
supported with range
$R$, that is, $C(x)=0$ for $|x|\geq R$.
\end{cond}

The following proposition shows that any kernel satisfying condition
$\mathcal{K}(\rho)
$ can be arbitrarily approximated by kernels verifying $ \mathcal
{K}_c(\rho,r)$
for $r$ large enough. We define the function $h$ by
\begin{equation}
\label{hfamily} h(x)= \exp \biggl(\frac{1}{|x|^2-1} \biggr) \mathbf{1}_{\lbrace|x| <
1 \rbrace
}
\qquad\forall x\in\mathbb{R}^d.
\end{equation}
For a function $f\in L^2(\mathbb{R}^d)$, put $\| f \|= \sqrt{\int
|f(t)|^2\,dt}$
and denote $[ f \ast f]$ the self-convolution product of $f$.

\begin{proposition}\label{convergencekernelcompactsupporttokernelverifyingK}
Let $C$ be a kernel verifying condition $\mathcal{K}(\rho)$ and $h$
be defined
by~(\ref{hfamily}). Then, for all $r>0$, the function $C_r$ defined by
\begin{equation}
\label{Csmoothcompactfamily} C_r(x)= \frac{1}{\|h \|^2} [ h \ast h] \biggl(
\frac{2x}{r} \biggr) C (x ) \qquad\forall x\in\mathbb{R}^d,
\end{equation}
verifies $\mathcal{K}_c(\rho,r)$. Moreover, we have the convergence
\begin{equation}
\label{convergenceCcompactsmooth} \lim_{r\rightarrow+\infty} C_r = C,
\end{equation}
uniformly on all compact sets.
\end{proposition}

In particular, by taking $C=C_B$ in~Proposition~\ref
{convergencekernelcompactsupporttokernelverifyingK}, it is always
possible to find a
kernel $C_r$ verifying $ \mathcal{K}_c(\rho,r)$ that yields a repulsiveness
(local or global) as close as we wish to the repulsiveness of $C_B$,
provided that $r$ is large enough. However, given a maximal range of
interaction $R$, it is clear that the maximal repulsiveness implied by
kernels verifying $\mathcal{K}_c(\rho,R)$ cannot reach the one of $C_B$,
since the support of $C_B$ is unbounded and $\operatorname{DPP}(C_B)$
is the unique
most repulsive DPP according to Corollary~\ref{mostrepulsiveDPP}. In
the following, we study the DPP's repulsiveness for a given range
$R>0$.

\begin{figure}

\includegraphics{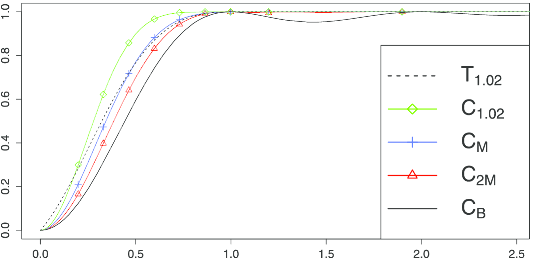}

\caption{In dimension $d=1$, comparison between the p.c.f. of
$\operatorname{DPP}(T_{1.02})$, $\operatorname{DPP}(C_{B})$ and
$\operatorname{DPP}(C_{R})$ for $R=1.02, M, 2M$.}
\label{figurecomparaisonturanCM}
\end{figure}

In comparison with condition $\mathcal{K}(\rho)$, the assumption that
$C$ is compactly
supported in condition $\mathcal{K}_c(\rho,R)$ makes the optimization problems
related to Definitions~\ref{globallyrepulsive}--\ref
{locallyrepulsive} much more difficult to investigate.
As a negative result, we know very little about the most globally
repulsive DPP, in the sense of Definition~\ref{globallyrepulsive},
under condition $\mathcal{K}_c(\rho,R)$. From relation~(\ref
{DPPpcf}), this is
equivalent to find a kernel $C$ with maximal $L^2$-norm under the
constraint that $C$ verifies $\mathcal{K}_c(\rho,R)$. Without
the constraint $\mathcal{F}(C)\leq1$, this problem is known as the
square-integral Tur\'an problem with range $R$ see, e.g.,
\cite{kolountzakis2003problem}. For this less constrained problem, it is
known that a solution exists, but no explicit formula is available,
cf.~\cite{Domar}.
For $d=1$, it has been proved that the solution is unique and there
exists an algorithm to approximate it; see~\cite{garsia1969}. In this
case, numerical approximations show that the solution with range $R$
verifies condition $\mathcal{K}_c(\rho,R)$ only if $R\leq1.02/\rho
$. This gives
the most globally repulsive DPP verifying $\mathcal{K}_c(\rho,R)$ in dimension
$d=1$, when $R\leq1.02/\rho$, albeit without explicit formula. Its
p.c.f.
is represented in Figure~\ref{figurecomparaisonturanCM}. For other
values of $R$, or in dimension $d\geq2$, no results are available, to
the best of our knowledge.

Let us now turn to the investigation of the most locally repulsive DPP,
in the sense of Definition~\ref{locallyrepulsive}, under condition
$\mathcal{K}_c
(\rho,R)$. Recall that without the compactly supported constraint of
the kernel, we showed in Section~\ref{sectionDPPregularitypcf} that the
most locally repulsive DPP, namely $\operatorname{DPP}(C_B)$, is also
(one of) the
most globally repulsive DPP.

For $\nu>0$, we denote by $j_\nu$ the first positive zero of the Bessel
function $J_\nu$ and by $J_\nu'$ the derivative of $J_\nu$. We refer to
\cite{stegun} for a survey about Bessel functions and their zeros.
Further, define the constant $M>0$ by
\[
M^d = \frac{2^{d-2} j^2_{({d-2})/{2}} \Gamma ({d}/{2}  )
}{\rho\pi^{{d}/{2}}}.
\]
We have $M\rho= \pi^2/8\approx1.234$ when $d=1$, $M \rho^{1/2} =
j_0/\pi^{1/2} \approx1.357$ when $d=2$ and $ M \rho^{1/3}= \pi^{1/3}
\approx1.465$ when $d=3$.

\begin{proposition}\label{compactDPP}
If $R\leq M$, then in the sense of Definition~\ref{locallyrepulsive},
there exists an unique isotropic kernel $C_R$ such that $\operatorname
{DPP}(C_R)$ is
the most locally repulsive DPP among all DPPs with kernel verifying
$\mathcal{K}_c(\rho,R)$. It is given by $C_R=u\ast u$ where
\begin{equation}
\label{u1} u(x)=\kappa \frac{J_{({d-2})/{2}}
(2j_{({d-2})/{2}}{|x|}/{R} )}{|x|^{{(d-2)}/{2}}} \mathbf {1}_{ \lbrace|x| <{R}/{2}
 \rbrace},
\end{equation}
with $\kappa^2=\frac{4 \Gamma (d/2 )}{ \rho\pi^{d/2}
R^2 }
(J'_{({d-2})/{2}}(j_{({d-2})/{2}}) )^{-2}$.
\end{proposition}

In this proposition $C_R$ is only given as a convolution product.
Nonetheless, an explicit expression is known in dimension $d=1$ and
$d=3$; see \cite{ehm2004convolution}. On the other hand, the Fourier
transform is known in any dimension since $\mathcal{F}(C_R)=\mathcal
{F}(u)^2$. We get
from the proof in Section~\ref{proofcompactDPP}, for all $x\in\mathbb{R}^d$,
\begin{equation}
\label{CR} \mathcal{F}(C_R) (x)=\rho\pi^{d/2}
R^d j^2_{({d-2})/{2}} \Gamma \biggl(\frac
{d}{2}
\biggr) \biggl( \frac{ J_{({d-2})/{2}}(\pi R|x|)} {  ( \pi
R|x|  )^{({d-2})/{2}}  ( j^2_{({d-2})/{2}}-  ( \pi
R|x|  )^2  )} \biggr)^2.
\end{equation}

If $R\geq M$, we have not been able to obtain a closed form expression
of the most locally repulsive stationary DPP. However, under some extra
regularity assumptions, we can state the following general result about
its existence and the form of the solution.

\renewcommand{\thecond}{$\mathcal{M}(\rho,R)$}
\begin{cond}
A function $u$ is said to verify condition $\mathcal{M}(\rho,R)$ if
$u(x)=0$ for
$|x|>\frac{R}{2}$, $u$ is a radial function and $u\in L^2(\mathbb
{R}^d)$ with
$\| u\|^2 = \rho$.
\end{cond}

\begin{proposition}\label{compactDPP2}
For any $R>0$, there exists an isotropic kernel $C_R$ such that
$\operatorname{DPP}(C_R)$ is the most locally repulsive DPP among all
DPPs with
kernel $C$ verifying $\mathcal{K}_c(\rho,R)$. It can be expressed as
$C_R=u\ast
u$ where $u$ satisfies $\mathcal{M}(\rho,R)$. Furthermore, if we
assume that
$\sup_{x\in\mathbb{R}^d} \mathcal{F}(C)(x) = \mathcal{F}(C)(0)$
and $u$ is twice differentiable
on its support, then $u$ is of the form
\begin{equation}
\label{u2} u(x) = \biggl(\beta+ \gamma\frac{J_{({d-2})/{2}} ( |x|/\alpha
 )}{|x|^{({d-2})/{2}}} \biggr)
\mathbf{1}_{ \lbrace|x|
<{R}/{2} \rbrace},
\end{equation}
where $\alpha>0,\beta\geq0$ and $\gamma$ are three constants linked
by the conditions $\mathcal{M}(\rho,R)$ and $\int_{\mathbb{R}^d}
u(x)\, dx \leq1 $.
\end{proposition}

In the case $R\leq M$, this proposition is a consequence of
Proposition~\ref{compactDPP} where $\beta=0$, $ \alpha
=R/(2j_{({d-2})/{2}})$ and $\gamma=\kappa$. When $R>M$, it is an open
problem to
find an explicit expression of the kernel $C_R$ without any extra
regularity assumptions. Even in this case, (\ref{u2}) only gives the
form of the solution and the constants $\alpha$, $\beta$ and $\gamma$
are not explicitly known. In particular, the choice $\beta=0$ does not
lead to the most locally repulsive DPP when $R>M$, contrary to the case
$R\leq M$.
In fact, the condition $\mathcal{M}(\rho,R)$ allows us to express
$\beta$ and
$\gamma$ as functions of $\alpha$, $R$ and $\rho$, but then some
numerical approximation are needed to find the value of $\alpha$ in
(\ref{u2}), given $R$ and~$\rho$, such that $\operatorname
{DPP}(C_R)$ is the most
locally repulsive DPP. We detail these relations in Section~\ref{compactfamilysection}, where we start from (\ref{u2}) to suggest a
new parametric family of compactly supported kernels.

Contrary to what happens in the noncompactly supported case of
Section~\ref{sectionDPPregularitypcf}, the most locally repulsive DPP
is not the most globally repulsive DPP under $\mathcal{K}_c(\rho,R)$.
This is
easily checked in dimension $d=1$ when
$R\leq1.02/\rho$ implying $R\leq M$: In this case the most globally
repulsive DPP under $\mathcal{K}_c(\rho,R)$ is $\operatorname
{DPP}(T_R)$, where $T_R$ is the
solution of the square-integral Tur\'an problem with range $R$, and the
most locally repulsive DPP is $\operatorname{DPP}(C_R)$ where $C_R$ is
given by (\ref{u1}). However, according to the results of Section~\ref{sectionDPPregularitypcf} corresponding to $R=\infty$, we expect that
$\operatorname{DPP}(C_R)$ has a strong global repulsiveness even for
moderate values
of $R$. This is confirmed in Figure~\ref{figurecomparaisonturanCM},
that shows the p.c.f. of $\operatorname{DPP}(C_R)$ when $d=1$, $\rho=1$
and $R=1.02$,
$R=M\approx1.234$ and $R=2M$, where in this case we take $C_R=u\ast u$
with $u$ given by (\ref{u2}) and the constants are obtained by
numerical approximations. The p.c.f.s of $\operatorname{DPP}(T_{1.02})$
and $\operatorname{DPP}(C_B)$
are added for sake of comparison.
Considering the behavior of the p.c.f. near the origin, we note that even
if $\operatorname{DPP}(T_{1.02})$ is the most globally repulsive DPP
under $\mathcal{K}_c(\rho
,R)$ when $R\leq1.02/\rho$, its local repulsiveness is not very
strong. On the other hand, $\operatorname{DPP}(C_R)$ seems to present
strong global
repulsiveness for the values of $R$ considered in the figure.

\section{Parametric families of DPP kernels}\label{DPPkernelsfamily}

A convenient parametric family of kernels $\{C_\theta\}_{\theta\in
\Theta
}$, where $\Theta\subset\mathbb{R}^q$ for some $q\geq1$, should ideally:
\begin{longlist}[(a)]
\item[(a)] provide a closed form expression for $C_\theta$, for any
$\theta$,
\item[(b)] provide a closed form expression for $\mathcal F(C_\theta)$,
for any $\theta$,
\item[(c)] be flexible enough to include a large range of DPPs, going
from the Poisson point process to $\operatorname{DPP}(C_B)$.
\end{longlist}
The second property above is needed to check the condition of existence
$\mathcal F(C_\theta)\leq1$, but it is also useful for some
approximations in practice. Indeed, the algorithm for simulating
$\operatorname{DPP}(C)$ on a compact set $S$, as presented in \cite
{hough2009zeros},
relies on the Mercer representation of $C$ on $S$, which is rarely
known in practice. In \cite{lavancierpublish}, this decomposition is
simply approximated by the Fourier series of~$C$ where, up to some
rescaling, the $k$th Fourier coefficients is replaced by
$\mathcal{F}
(C)(k)$. The same approximation is used to compute the likelihood. This
method has proved to be accurate in most cases, both from a practical
and a theoretical point of view, provided $\rho$ is not too small, and
to be computationally efficient; see~\cite{lavancierpublish}.

In addition to (a)--(c), we may also require that $C_\theta$ is
compactly supported with maximal range~$R$, following the motivation
explained in Section~\ref{sectionDPPregularitypcffiniterange}, in
which case the maximal possible repulsiveness is given by
$\operatorname{DPP}(C_R)$.
Or we may require that $\mathcal F(C_\theta)$ is compactly supported,
in which case the Fourier series mentioned in the previous paragraph
becomes a finite sum and no truncation is needed in practice. Note,
however, that $C_\theta$ and $\mathcal F(C_\theta)$ cannot both be
compactly supported.

Several standard parametric families of kernels are available,
including the well-known Whittle--Mat\'ern and the generalized Cauchy
covariance functions, where the condition $\mathcal F(C_\theta)\leq1$
implies some restriction on the parameter space; see~\cite
{lavancierpublish}. Although they encompass a closed form expression
for both $C_\theta$ and $\mathcal F(C_\theta)$, they are not flexible
enough to reach the repulsiveness of $\operatorname{DPP}(C_B)$.
Another family of
parametric kernels is considered in \cite{lavancierpublish}, namely
the power exponential spectral model, that contains as limiting cases
$C_B$ and the Poisson kernel (\ref{kernelpoisson}). For this reason,
this family is more flexible than the previous ones, but then only
$\mathcal F(C_\theta)$ is given and no closed expression is available
for $C_\theta$. For all these families, none of $C_\theta$ and
$\mathcal F(C_\theta)$ is compactly supported.

Below, we present alternative families of parametric kernels. The first
two ones, so-called Bessel-type and Laguerre--Gaussian families, fulfil
the three requirements (a)--(c) above and the Bessel-type family has the
additional property that the Fourier transform of the kernels is
compactly supported. Moreover, we introduce new families of compactly
supported kernels, inspired by Proposition~\ref
{convergencekernelcompactsupporttokernelverifyingK} and
Proposition~\ref{compactDPP2}.

\subsection{Bessel-type family}\label{besselfamily}

For all $\sigma\geq0, \alpha>0, \rho>0 $, we consider the
Bessel-type kernel
\begin{equation}
\label{bessel-type} C(x)= \rho 2^{({\sigma+d})/{2}}\Gamma \biggl(\frac{\sigma+d+2}{2}
\biggr)\frac{J_{({\sigma+d})/{2}} (2 |{x}/{\alpha
}|\sqrt{({\sigma+d})/{2}} )}{ (2|{x}/{\alpha}|\sqrt{({\sigma+d})/{2}} )^{({\sigma+d})/{2}}},\qquad x\in\mathbb{R}^d.
\end{equation}
This positive definite function first appears in \cite
{schoenberg1938metric}, where it is called the Poisson function. It has
been further studied in \cite{flyer2006exact} and \cite
{fornberg2006}, where it is called the Bessel-type function. For
obvious reasons, we prefer the second terminology when applied to point
processes. For any $x\in\mathbb R$, we denote by $x_+=\max(x,0)$ its
positive part.

\begin{proposition}\label{DPPB}
Let $C$ be given by (\ref{bessel-type}), then its Fourier transform is,
for all $x\in\mathbb{R}^d$,
\begin{equation}
\label{Fouriertransformgeneralizedjinc}
\mathcal{F}(C) (x)= \rho \frac{ (2\pi)^{{d}/{2}}\alpha^d
\Gamma(({\sigma+d+2})/{2}) }{(\sigma+d)^{{d}/{2}} \Gamma(({\sigma+2})/{2})} \biggl(1-
\frac{2 \pi^2 \alpha^2 |x|^2}{\sigma+d} \biggr)^{{\sigma}/{2}}_{+}
\end{equation}
and $\operatorname{DPP}(C)$ exists if and only if $\alpha\leq\alpha
_{\max}$ where
\[
\alpha_{\max}^d=\frac{(\sigma+d)^{{d}/{2}} \Gamma(({\sigma+2})/{2})}{\rho(2\pi)^{{d}/{2}}
\Gamma(({\sigma+d+2})/{2})}.
\]
In this case, $\operatorname{DPP}(C)$ defines a stationary and
isotropic DPP with
intensity $\rho$.
Moreover, if $\sigma=0$ and $\alpha= \alpha_{\max}$, then $C=C_{B}$
where $C_{B}$ is defined in~(\ref{expressionjinc}). In addition, for
any $\rho>0$ and $\alpha>0$, we have the convergence
\begin{equation}
\label{Bconvergence} \lim_{{\sigma\rightarrow+\infty}} C(x) = \rho e^{- ( {|x|}/{\alpha}  )^2},
\end{equation}
uniformly on all compact sets.
\end{proposition}

The Bessel-type family contains $C_B$ as a particular case and the
Poisson kernel as a limiting case, when $\alpha\to0$. Moreover,
$\mathcal{F}
(C)$ is compactly supported; see~(\ref
{Fouriertransformgeneralizedjinc}). The plots in Figure~\ref{figureoverview}(b)--(d) show some
realizations of this model when $\sigma=0$ and $\alpha=0.2, 0.4,
\alpha_{\max}$, respectively. The supplementary material~\cite
{bisciosupplementary} includes more simulations and shows the behavior
of the p.c.f. for different values of the parameters.

\subsection{Laguerre--Gaussian family}

Let us first recall the definition of the Laguerre polynomials. We
denote by $\mathbb{N}$ the set $\{0,1,2,\ldots\}$\vspace*{1pt} and by $\mathbb
{N}^*$ the set $\mathbb{N}
\setminus\{0\}$. For integers $0\leq k \leq m $ and numbers $\alpha$,
define ${m+\alpha\choose k}=\frac{(m+\alpha)\cdots(m+\alpha+1-k)}{k}$
if $k>0$ and ${m+\alpha\choose k}=1$ if $k=0$.

\begin{definition}
The Laguerre polynomials are defined for all $m\in\mathbb{N}$ and
$\alpha\in\mathbb{R}
$ by
\[
L_m^\alpha(x) = \sum_{k=0}^m
\pmatrix{m+\alpha\cr m-k} \frac
{(-x)^k}{k!}\qquad\forall x\in\mathbb{R}.
\]
\end{definition}

For all $m\in\mathbb{N}^{*}, \alpha>0, \rho>0 $ and $x\in\mathbb
{R}^d$, we
consider the Laguerre--Gaussian function
\begin{equation}
\label{laguerregausskernel} C (x) = \frac{\rho}{{m-1+{d}/{2} \choose m-1} } L_{m-1}^{{d}/{2}} \biggl(
\frac{1}{m}\biggl\llvert \frac{x}{\alpha}\biggr\rrvert ^2
\biggr)e^{-{1}/{m} |{x}/{\alpha}|^2}.
\end{equation}
This kernel already appears in the literature; see, for example, \cite
{fasshauer2007} for an application in approximation theory. The
following proposition summarizes the properties that are relevant for
its use as a DPP kernel.

\begin{proposition}\label{DPPlaguerre}
Let $C$ be given by (\ref{laguerregausskernel}), then its Fourier
transform is, for all $x\in\mathbb{R}^d$,
\begin{equation}
\label{Fourierlaguerrefamily} \mathcal{F} (C ) (x)= \frac{\rho}{{{m-1+{d}/{2}}\choose{m-1}}} \alpha ^d (m
\pi )^{{d}/{2}} e^{-m(\pi\alpha|x|)^2} \sum_{k=0}^{m-1}
\frac{(\pi\sqrt{m} |\alpha x|)^{2k}}{k!}
\end{equation}
and $\operatorname{DPP}(C)$ exists if and only if $\alpha\leq\alpha
_{\max}$ where
\[
\alpha_{\max}^d=\frac{{{m-1+{d}/{2}}\choose{m-1}}}{\rho(m\pi
)^{{d}/{2}}}.
\]
In this case, $\operatorname{DPP}(C)$ is stationary and isotropic with
intensity $\rho
$. Moreover, for any $\rho>0$ and $\alpha>0$, we have the convergence
\begin{equation}
\label{LGconvergence1} \lim_{m\rightarrow+\infty} C(x) = \rho\Gamma \biggl(
\frac
{d}{2}+1 \biggr) \frac{J_{{d}/{2}} (2|{x}/{\alpha}| )}{ |{x}/{\alpha}|^{{d}/{2}}}
\end{equation}
uniformly on all compact sets. In particular, for $\alpha=\alpha
_{\max}$,
\begin{equation}
\label{LGconvergence2} \lim_{m\rightarrow+\infty} C(x) = C_B(x)
\end{equation}
uniformly on all compact sets and where $C_B$ is defined in~(\ref
{expressionjinc}).
\end{proposition}

This family of kernels contains the Gaussian kernel, being the
particular case $m=1$, and includes as limiting cases the Poisson
kernel~(\ref{kernelpoisson}) (when $\alpha\to0$) and $C_B$, in view
of (\ref{LGconvergence2}). Some illustrations of this model are
provided in the supplementary material~\cite{bisciosupplementary},
including graphical representations of the p.c.f. and some realizations.

\subsection{Families of compactly supported kernels}\label
{compactfamilysection}

As suggested by Proposition~\ref
{convergencekernelcompactsupporttokernelverifyingK}, we can consider
the following family of compactly
supported kernels, parameterized by the range $R>0$,
\begin{equation}
\label{parametricsmoothcompactfamily} C_1(x)= \frac{1}{\|h \|^2} [ h \ast h] \biggl(
\frac{2x}{R} \biggr) C_B (x )\qquad\forall x\in
\mathbb{R}^d,
\end{equation}
where $h$ is given by~(\ref{hfamily}). The Poisson kernel (\ref
{kernelpoisson}) and $C_B$ are two limiting cases, when respectively
$R\to0$ and $R\to+\infty$. However, this family of kernels has
several drawbacks: No closed form expression is available for $C_1$,
nor for $\mathcal{F}(C_1)$. Moreover, when the range $R$ is fixed,
$\operatorname{DPP}(C_1)$ is
not the most repulsive DPP; see Proposition~\ref{compactDPP2} and the
graphical representations in the supplementary material~\cite
{bisciosupplementary}. This is the reason why we turn to another family
of compactly supported kernels.

Following Proposition~\ref{compactDPP2}, we introduce a new family of
compactly supported kernels with range $R$, given as a convolution
product of functions as in (\ref{u2}). Specifically, let $R>0$, $\rho
>0$ and $\alpha>0$ such that $R/(2\alpha)$ is not a zero of the Bessel
function $J_{({d-2})/{2}}$ and consider the kernel $C_2=u\ast u$ with
\begin{equation}
\label{familyu} u(x)= \sqrt{\rho} \beta(R,\alpha) \biggl( 1-\frac{R^{{d}/{2}-1}}{
2^{{d}/{2}-1} J_{{d}/{2}-1}({R}/({2\alpha}))}
\frac{J_{{d}/{2}-1}(\llvert {x}/{\alpha}\rrvert )}{|x|^{{d}/{2}-1}} \biggr) \mathbf{1} _{  \lbrace|x|\leq{R}/{2}  \rbrace},
\end{equation}
where
\begin{eqnarray*}
\beta(R,\alpha)
& =&  \biggl[ \frac{R^{d-1}\pi^{d/2}}{ 2^{d-1}
\Gamma
({d}/{2} ) } \biggl( \frac{R}{d}-4\alpha\frac{J_{{d}/{2}}({R}/({2\alpha}))}{J_{{d}/{2}-1} ({R}/({2\alpha}) )}\\
&&{} +
\frac{R}{2} \biggl(1-\frac{J_{{d}/{2}-2}({R}/({2\alpha}))J_{{d}/{2}}({R}/({2\alpha}))}
{J_{{d}/{2}-1}^2 ({R}/({2\alpha}) )} \biggr)
\biggr) \biggr]^{-{1}/{2}}.
\end{eqnarray*}

\begin{proposition}\label{DPPcompact}
Let $C_2=u\ast u$ where $u$ is given by (\ref{familyu}), then its
Fourier transform is $\mathcal F(u)^2$ where for all $x\in\mathbb{R}^d$
\begin{eqnarray*}
&& \!\!\!\mathcal{F}(u) (x)\\
 &&\!\!\!\quad= \sqrt\rho \beta(R,\alpha) \biggl(\frac{R}{2|x|}
\biggr)^{{d}/{2}-1} \biggl(\frac{R}{2|x|} J_{{d}/{2}} \bigl(\pi R |x|\bigr)
+ \frac{\pi}{J_{({d-2})/{2}}({R}/({2\alpha}))}\\
&&\!\!\!\qquad {}\times\frac{ R\alpha J_{({d-2})/{2}}'({R}/({2\alpha}))
J_{({d-2})/{2}}(\pi R|x|)-2 \pi R\alpha^2 J_{({d-2})/{2}}({R}/({2\alpha}))|x| J_{({d-2})/{2}}'
(\pi R|x|)}{1-4\pi^2 |\alpha x|^2} \biggr).
\end{eqnarray*}
Moreover, $\operatorname{DPP}(C_2)$ exists if and only if $\alpha$ is
such that
$|\mathcal F(u)|\leq1$. In this case, $\operatorname{DPP}(C_2)$
defines a stationary
and isotropic $R$-dependent DPP with intensity $\rho$.
\end{proposition}

The choice of $u$ in (\ref{familyu}) comes from (\ref{u2}) where
$\gamma$ has been chosen such that $u$ is continuous at $|x|=R/2$ and
where $\beta$ is deduced from the relation $C_2(0)=\|u\|^2=\rho$. Given
$\rho$ and $R$, the remaining free parameter in this parametric family
becomes $\alpha$.
The restriction that $R/(2\alpha)$ must not be a zero of $J_{({d-2})/{2}}$ can be alleviated by setting in these cases $\beta=0$ in
(\ref{u2}) and choose $\gamma$ so that $C_2(0)=\rho$. Then the most
locally repulsive DPP (\ref{u1}) when $R\leq M$ would be part of the
parametric family. However, these kernels can be arbitrarily
approximated by some kernel given by (\ref{familyu}) for some value of
$\alpha$, so we do not include these particular values of $\alpha$ in
the family above.

The condition $|\mathcal F(u)|\leq1$ on $\alpha$, given $R$ and $\rho
$, must be checked numerically. In most cases, the maximal value of
$\mathcal F(u)$ holds at the origin and we simply have to check whether
\mbox{$|\mathcal F(u)(0)|\leq1$}. No theoretical results are available to
claim the existence of an admissible $\alpha$, but from our experience,
there seems to exist an infinity of admissible $\alpha$ for any $R$ and
$\rho$. Moreover, while the most locally repulsive DPP when $R\leq M$
is known and corresponds to (\ref{u1}), the most repulsive DPP when
$R>M$ in the above parametric family seems to correspond to the maximal
value of $\alpha$ such that $|\mathcal F(u)|\leq1$, denoted $\alpha
_{\max}$.

The parametric family given by $C_2$ is mainly of interest since it
covers a large range of repulsive DPPs while the kernels are compactly
supported. Moreover, the closed form expression of $\mathcal F(C_2)$ is
available and this family contains the most locally repulsive DPP with
range $R$, in view of Proposition~\ref{compactDPP2}, at least when
$R\leq M$. Some illustrations are provided in the supplementary
material~\cite{bisciosupplementary}.

\section{Proofs}\label{proofs}

\subsection{Proof of Proposition~\texorpdfstring{\protect\ref{maxlocrep}}{3.2}}

As the kernel $C_{B}$ verifies condition $\mathcal{K}(\rho)$, it
defines a $\mathrm{DPP}$ with
intensity $\rho$ and its associated p.c.f. $g_B$ given by~(\ref{DPPpcf})
vanishes at $0$.
By the analytic definition of Bessel functions; see \cite{stegun},
relation~(9.1.10)],
\begin{eqnarray*}
&& C_{B}(x)= \frac{\sqrt{\rho\Gamma({d}/{2}+1)}}{\pi^{d/4}} \sum_{n=0}^{+\infty}
\frac{(-1)^n  ( \sqrt{\pi}\Gamma({d}/{2}+1)^{{1}/{d}}\rho^{{1}/{d}}
)^{2n}}{2^{2n}n!\Gamma
 (n+1+d/2 )} |x|^{2n}.
\end{eqnarray*}
Thus, $C_{B}$ is twice differentiable at $0$ and by (\ref{DPPpcf}), so
is $g_{B}$.
By Definition~\ref{locallyrepulsive}, any DPP having a p.c.f. $g$ that
does not vanish at 0 or is not twice differentiable at $0$ is less
locally repulsive than $\operatorname{DPP}(C_{B})$. Consequently, we
assume in the
following of the proof that $g(0)=0$ and $g$ is twice differentiable at $0$.
The problem therefore reduces to minimize $\Delta g(0)$ under the
constraint that $g$ is the p.c.f. of a DPP with kernel $C$ verifying
condition $\mathcal{K}(\rho)$.

According to condition $\mathcal{K}(\rho)$, the Fourier transform of
the kernel $C$ is
well defined and belongs to $L^1(\mathbb{R}^d)$, as noticed below
(\ref{definitionfouriertransform}).
Therefore, we can define the function $f=\frac{\mathcal{F}(C)}{\|
\mathcal{F}(C)\|_1}$
where $\|\mathcal{F}(C)\|_1=\int_{\mathbb{R}^d} |\mathcal
{F}(C)(x)|\,dx$ and consider it as a
density function of a random variable $X=(X_1,\ldots,X_d)\in\mathbb
{R}^d $.
Denote by $\widehat{f}(t) = \mathbb{E} (e^{i t\cdot X} )$ the
characteristic function of $X$. We have
\begin{eqnarray}
\label{lienfourierprobavsanalyse} \widehat{f} (t) &=& \frac{C ({t}/({2\pi}) )}{\|
\mathcal{F}(C) \|
_1}\qquad\forall t\in
\mathbb{R}^d.
\end{eqnarray}
Thus, $\widehat{f}$ is twice differentiable at $0$, so by the usual
properties of the characteristic function (see~\cite
{shiryaev1995probability}), $X$ has finite second-order moments and
\begin{eqnarray}
\label{expressioncoordonnemoments} E\bigl(X_i^2\bigr)&=&-
\frac{\partial^2 \widehat{f}}{\partial x_i^2}(0)+ \biggl(\frac
{\partial\widehat{f}}{\partial x_i}(0) \biggr)^2,\qquad
i=1, \ldots, d.
\end{eqnarray}

On the other hand, as already noticed in Section~\ref{sectionintroduction}, $\nabla g(0) =0$ and so $\frac{\partial
C}{\partial
x_i}(0)=0$ for $i=1,\dots,d$.
By differentiating both sides of (\ref{lienfourierprobavsanalyse}),
\begin{eqnarray}
\label{partialdef} \frac{\partial\widehat{f}}{\partial x_i}(0) &=&\frac{1}{2\pi\|
\mathcal{F}(C)\|
_1} \frac{\partial C}{\partial x_i}(0)=0,
\qquad i=1, \ldots, d
\end{eqnarray}
and
\begin{eqnarray}
\label{partial2def} \frac{\partial^2 \widehat{f}}{\partial x_i^2}(0) &=&\frac{1}{4\pi
^2 \|\mathcal{F}
(C)\|_1} \frac{\partial^2 C}{\partial x_i^2}(0),
\qquad i=1, \ldots, d.
\end{eqnarray}
Then, by (\ref{expressioncoordonnemoments})--(\ref{partial2def}),
\begin{eqnarray*}
E\bigl(|X|^2\bigr) &=& E \Biggl(\sum_{i=1}^d
X_i^2 \Biggr) = - \Delta\widehat{f}(0) = -
\frac{1}{4\pi^2 \|\mathcal{F}(C)\|_1} \Delta C(0).
\end{eqnarray*}
Moreover,
\begin{eqnarray*}
E\bigl(|X|^2\bigr) &=& \int_{\mathbb{R}^d}
|x|^2 f(x) \,dx = \int_{\mathbb
{R}^d} |x|^2
\frac{\mathcal{F}
(C)}{\|\mathcal{F}(C)\|_1}(x)\, dx.
\end{eqnarray*}
Hence,
\begin{eqnarray}
\label{laplaciendeC} \Delta C(0)&=& -4\pi^2\int_{\mathbb{R}^d}
|x|^2 \mathcal{F}(C) (x) \,dx.
\end{eqnarray}
By (\ref{DPPpcf}) and since $\nabla C(0)=0$,
\begin{eqnarray}
\Delta g(0)&=& \Delta \biggl(1-\frac{C^2}{\rho^2} \biggr) (0) = -
\frac
{1}{\rho
^2} \Biggl( \sum_{i=1}^{d}
2C(0) \frac{\partial^2 C}{\partial x_i^2}(0)
+ 2 \biggl(\frac{\partial C}{\partial x_i}(0) \biggr)^2 \Biggr)
\nonumber
\\[-8pt]
\label{laplaciengenfonctiondeC}
\\[-8pt]
\nonumber
&=& -\frac{2}{\rho} \sum_{i=1}^{d}
\frac{\partial^2 C}{\partial x_i^2}(0) = -\frac{2}{\rho} \Delta C(0).
\end{eqnarray}
Finally, we deduce from (\ref{laplaciendeC}) and (\ref
{laplaciengenfonctiondeC}) that
\[
\Delta g(0) = \frac{8\pi^2}{\rho} \int_{\mathbb{R}^d} |x|^2
\mathcal{F}(C) (x)\,dx.
\]
Thus, the two following optimization problems are equivalent.

\begin{problem}
Minimizing $\Delta g(0)$ under the constraint
that $g$ is the p.c.f. of a DPP with kernel $C$ satisfying condition
$\mathcal{K}(\rho)$.
\end{problem}

\begin{problem}
Minimizing $\int_{\mathbb{R}} |x|^2 \mathcal{F}(C)(x)\,dx$
under the constraint that $C$ is a kernel which is twice differentiable
at $0$ and verifies the condition $\mathcal{K}(\rho)$.
\end{problem}

The latter optimization problem is a special case of \cite
{liebanalysis2001}, Theorem~1.14, named bathtub principle, which gives the
unique solution $\mathcal{F}(C) = \mathbf{1}_{ \lbrace|\cdot
|^d \leq\rho\tau^d
 \rbrace}$ in agreement with (\ref{expressionjinc}). This
completes the proof.

\subsection{Proof of Proposition~\texorpdfstring{\protect\ref
{convergencekernelcompactsupporttokernelverifyingK}}{4.1}}

Notice that $h$ is symmetric, real-valued, infinitely differentiable
and verifies $h(x)=0$ for $x\geq1$; see \cite
{sasvari2013multivariate}, Section~3.2. Thus, $\|h\|$ is finite and $\|
h\|\neq
0$, so $C_r$ is well defined.

Since $h \ast h(0) = \| h\|^2$, we have $C_r(0)=\rho$. By product
convolution properties, $h\ast h$ is symmetric, real-valued, infinitely
differentiable and compactly supported with range $2$. Thus, by~(\ref
{Csmoothcompactfamily}), $C_r$ is symmetric, real-valued, infinitely
differentiable and compactly supported with range $r$.
Then $C_r$ belongs to $L^1(\mathbb{R}^d) \cap L^2(\mathbb{R}^d) $. In
particular, $\mathcal{F}
(C_r)$ is well defined pointwise. By well-known properties of the
Fourier transform, for all $x\in\mathbb{R}^d$,
\begin{equation}
\label{fourierCcompactsmooth} \mathcal{F}(C_r) (x) =\frac{r^d}{2^d\|h\|^2} \biggl[
\mathcal {F}(h)^2 \biggl(\frac{r }{2} \cdot \biggr) \ast
\mathcal{F}(C) ( \cdot) \biggr](x).
\end{equation}
Since $h$ is symmetric, $\mathcal{F}(h)$ is real valued, so $\mathcal
{F}(h)^2\geq0$.
Thus, as $\mathcal{F}(C) \geq0$ by condition $\mathcal{K}(\rho)$,
we have $\mathcal{F}(C_r) \geq0$.
Further, since $0\leq\mathcal{F}(C)\leq1$,
\begin{equation}
\label{majorationCcompactsmooth} \frac{r^d}{2^d\|h\|^2} \int_{\mathbb{R}^d}
\mathcal{F}(h)^2 \biggl(\frac{r t}{2} \biggr) \mathcal{F} (C)
(x-t) \,dt \leq\frac{r^d}{2^d\|h\|^2} \int_{\mathbb{R}^d}
\mathcal{F}(h)^2 \biggl(\frac
{r t}{2} \biggr) \,dt.
\end{equation}
By the substitution $u=rt/2$ and Parseval's equality, the right-hand
side of (\ref{majorationCcompactsmooth}) equals $1$. Finally, (\ref
{fourierCcompactsmooth}) and (\ref{majorationCcompactsmooth}) give
$ \mathcal{F}(C_r) \leq1$, that is, $0\leq\mathcal{F}(C_r) \leq1$.

It remains to show the convergence result~(\ref
{convergenceCcompactsmooth}), which reduces to prove
that $ \frac{1}{\|h \|^2} [ h\ast h] ( \frac{2}{r} \cdot
) $
tends to $1$ uniformly on all compact set when $r\to\infty$.
This follows from $h \ast h (0) = \|h\|^2$ and the uniform continuity
of $h\ast h$ on every compact set.

\subsection{Proof of Proposition~\texorpdfstring{\protect\ref{compactDPP}}{4.2}} \label
{proofcompactDPP}

The proof is based on a theorem from
Ehm, Gneiting and Richards \cite{ehm2004convolution}
recalled below with only slight changes in the presentation.

\begin{definition}\label{radialization}
Let $H$ denote the normalized Haar measure on the group $\mathrm{SO}(d)$ of
rotations in $\mathbb{R}^d$ and let $C$ be a kernel verifying
condition $\mathcal{K}_c
(\rho,R)$. The radialization of the kernel $C$ is the kernel $\operatorname{rad}(C)$ defined by
\begin{eqnarray*}
\operatorname{rad}(C) (x) &=& \int_{\mathrm{SO}(d)} C\bigl( j(x) \bigr)
H(dj).
\end{eqnarray*}
\end{definition}

Note that for any isotropic kernel $C$, $C=\operatorname{rad}(C)$. We
say that
$C_1=C_2$ up to a radialization if $C_1$ and $C_2$ are kernels
verifying condition $\mathcal{K}_c(\rho,R)$ and $\operatorname
{rad}(C_1)=\operatorname{rad}(C_2)$.

Define $\gamma_d >0$ by $\gamma_d^2=\frac{4j^{d-2}_{(d-2)/2}}{\pi
^{{d}/{2}}\Gamma ({d}/{2} )J^2_{{d}/{2}}(j_{(d-2)/2})}$
and set $c_d=\frac{4j^2_{({d-2})/{2}}}{4^d\pi^{{d}/{2}}\Gamma
({d}/{2} )}$ where $j_{(d-2)/2}$ is introduced before
Proposition~\ref{compactDPP}.

\begin{theorem}[(Ehm, Gneiting and Richards \cite
{ehm2004convolution})] \label{ehmthm}
Let $\Psi$ be a twice differentiable characteristic function of a
probability density $f$ on $\mathbb{R}^d$ and suppose that $\Psi
(x)=0$ for
$|x|\geq1$. Then
\[
-\Delta\Psi(0)=\int|x|^2 f(x)\,dx \geq4 j_{(d-2)/2}^2
\]
with equality if and only if, up to a radialization, $\Psi= \omega_d
\ast\omega_d$, where
\[
\omega_d(x)= %
\cases{ \gamma_d
\frac{\Gamma ({d}/{2}
)}{j_{(d-2)/2}^{(d-2)/2}} \frac{J_{({d-2})/{2}} (2j_{({d-2})/{2}}
|x|)}{ |x|^{({d-2})/{2}}}, & \quad$\mbox{if } |x|\leq
\frac
{1}{2}$,\vspace*{3pt}
\cr
0, & \quad\mbox{otherwise}.} %
\]
The corresponding minimum variance density is
\[
f(x) = c_d \Gamma \biggl(\frac{d}{2} \biggr)^2
\biggl( \frac{2^{({d-2})/{2}} J_{({d-2})/{2}}({|x|}/{2})}{ |{x}/{2}|^{({d-2})/{2}}  ( j^2_{(d-2)/2}- ({|x|}/{2} )^2  )} \biggr)^2.
\]
\end{theorem}
According to Definition~\ref{locallyrepulsive} and by the same
arguments as in the proof of Proposition~\ref{maxlocrep} and~(\ref
{laplaciengenfonctiondeC}), we seek a kernel $C$ which is twice
differentiable at $0$ such that $\Delta C(0)$ is maximal among all
kernels verifying condition $\mathcal{K}_c(\rho,R)$.

In a first step, we exhibit a candidate for the solution to this
optimization problem and in a second step we check that it verifies all
required conditions.

\begin{longlist}[\textit{Step} 1.]
\item[\textit{Step} 1.]
We say that a function $C$ verifies $\widetilde{\mathcal{K}_c}(\rho
,R)$ if it
verifies $\mathcal{K}_c(\rho,R)$ without necessarily verifying
$\mathcal{F}(C)\leq1$.
Notice that a function $C$ verifies $\widetilde{\mathcal{K}_c}(\rho
,R)$ if and
only if the function
\begin{equation}
\label{lienPsietC} \Psi(x) = \frac{C(Rx)}{\rho}, \qquad x\in\mathbb{R}^d,
\end{equation}
verifies $\widetilde{\mathcal{K}_c}(1,1)$. Therefore, we have a one-to-one
correspondence between $\widetilde{\mathcal{K}_c}(\rho,R)$ and
$\widetilde{\mathcal{K}_c}(1,1)$.

On the other hand, if a function $\Psi$ verifies condition $\widetilde
{\mathcal{K}_c}(1,1)$, it is by Bochner's theorem the characteristic
function of
a random variable $X$. Moreover, the function $\Psi$ is continuous and
compactly supported, so it is in $L^1(\mathbb{R}^d)$ and the random variable
$X$ has a density $f$; see \cite{shiryaev1995probability}.
Thus, by Theorem~\ref{ehmthm}, any function $\Psi$ twice differentiable
at $0$ and verifying condition $\widetilde{\mathcal{K}_c}(1,1)$ satisfies
\begin{equation}
\label{majorationcasO1} \Delta\Psi(0) \leq-4 j_{(d-2)/2}^2.
\end{equation}
By differentiating both sides of (\ref{lienPsietC}), we have
\begin{equation}
\label{relationlaplacienpsietC} \Delta\Psi(0) = \frac{R^2}{\rho}\Delta C(0).
\end{equation}
Thus, by (\ref{majorationcasO1}) and (\ref{relationlaplacienpsietC}),
for any kernel $C$ which is twice differentiable at $0$ and verifies
$\widetilde{\mathcal{K}_c}(\rho,R)$,
\begin{equation}
\label{majorationcasOR} \Delta C(0)=\frac{\rho\Delta\Psi(0)}{R^2} \leq-\frac{4\rho
j_{(d-2)/2}^2 }{R^2}.
\end{equation}
By Theorem~\ref{ehmthm}, the equality in~(\ref{majorationcasOR})
holds if and only if $\Psi= \omega_d \ast\omega_d$ and we name $C_R$
the corresponding kernel $C$ given by (\ref{lienPsietC}).

\item[\textit{Step} 2.] The kernel $C_R$ is the candidate to our
optimization problem, however, it remains to prove that it verifies
condition $\mathcal{K}_c(\rho,R)$.
We have seen in step~1 that $C_R$ verifies $\widetilde
{\mathcal{K}_c
}(\rho,R)$ and is twice differentiable at $0$. We must show that
$\mathcal{F}
(C_R) \leq1$.
By Theorem~\ref{ehmthm}, the function $\Psi= \omega_d \ast\omega
_d$ is
the characteristic function of a probability density $f$. Thus, for all
$x\in\mathbb R^d$,
\begin{equation}
\label{transformeefourierdePsi} \mathcal{F} (\Psi ) (x)=(2\pi)^d f (2\pi x ) = (2
\pi)^d c_d \Gamma \biggl(\frac{d}{2}
\biggr)^2 \biggl( \frac{2^{({d-2})/{2}}
J_{({d-2})/{2}}(|\pi x|)}{ |\pi x|^{({d-2})/{2}}  (
j^2_{(d-2)/2}-
(|\pi x| )^2  )} \biggr)^2.
\end{equation}
By (\ref{lienPsietC}) and the Fourier transform dilatation we
thereby obtain (\ref{CR}).

Moreover, the Bessel functions are nonnegative up to their first
nonnegative zero so $\omega_d\geq0$, which implies that $\Psi\geq0$.
Hence, by (\ref{transformeefourierdePsi}),
\begin{equation}
\label{majorationfourierPsi} \mathcal{F}(\Psi) (x)=\biggl\llvert \int_{\mathbb{R}^d}
\Psi(t) e^{2i\pi
x\cdot t}\,dt\biggr\rrvert \leq \int_{\mathbb{R}^d}
\Psi(t)\,dt=\mathcal{F}(\Psi) (0)=\frac{2^d \pi
^d c_d}{j^4_{({d-2})/{2}}}.
\end{equation}
Thus, by (\ref{lienPsietC}) and the Fourier transform dilatation,
\begin{equation}
\mathcal{F}(C_R) (x)\leq\mathcal{F}(C_R) (0)=
\frac{2^d R^d \rho\pi
^d c_d}{j^4_{({d-2})/{2}}}=\frac{R^d}{M^d}.
\end{equation}
Since by hypothesis $R \leq M$, we have $\mathcal{F}(C_R)\leq1$.
\end{longlist}

\subsection{Proof of Proposition~\texorpdfstring{\protect\ref{compactDPP2}}{4.3}} \label
{proofcompactDPP2}

According to Definition~\ref{locallyrepulsive} and by the same
arguments as in the proof of Proposition~\ref{maxlocrep} and~(\ref
{laplaciengenfonctiondeC}), we seek a kernel $C$ which is twice
differentiable at $0$ such that $\Delta C(0)$ is maximal among all
kernels verifying condition $\mathcal{K}_c(\rho,R)$. By~(\ref{laplaciendeC}),
this is equivalent to solve the following Problem~\ref{A}.

\renewcommand{\theproblem}{A}
\begin{problem}\label{A}
Minimize $\int_{\mathbb{R}^d} |x|^2 \mathcal{F}(C)(x)\,dx$ under the
constraints that $C$
is twice differentiable at~$0$ and verifies $\mathcal{K}_c(\rho,R)$.
\end{problem}

The proof of Proposition~\ref{compactDPP2} is based on the following
three lemmas. In the first lemma, the gradient $\nabla u$ has to be
considered in the sense of distribution when $u\in L^2(\mathbb{R}^d)$
is not
differentiable.

\begin{lemma}\label{minimisationproblem}
A kernel $C_R$ is solution to Problem~\ref{A} if and only if there exists a
function $u$ such that, up to a radialization, $C_R=u\ast u$ where $u$
minimizes $\int_{\mathbb{R}^d} |\nabla u (x) |^2 \,dx$ among all
functions $u$
verifying $\mathcal{M}(\rho,R)$ and $\mathcal{F}(u)^2 \leq1 $.
\end{lemma}

The existence statement in Proposition~\ref{compactDPP2} is given by
the following lemma.

\begin{lemma}\label{existenceminimsobolev}
There exists a solution to Problem~\ref{A}.
\end{lemma}

By Lemma~\ref{minimisationproblem}, $C_R=u \ast u$ where $u$ is the
solution of the given optimization problem.
Then, under the additional constraint $\sup_{x\in\mathbb{R}^d}
\mathcal{F}(C)(x) = \mathcal{F}
(C)(0)$, we have $\sup_{x\in\mathbb{R}^d} ( \mathcal{F}(u)(x)
 )^2 =  (\mathcal{F}
(u)(0) )^2$. Since $ \mathcal{F}(u)^2(0) =  ( \int_{\mathbb{R}^d} u(t) \,dt
)^2$, the constraint $\mathcal{F}(u)^2 \leq1 $ in Lemma~\ref
{minimisationproblem} becomes $ ( \int_{\mathbb{R}^d} u(t) \,dt
 )^2\leq1$.
Notice that $-u$ is also a solution of the optimization problem. Thus,
we can assume without loss of generality that $\int_{\mathbb{R}^d}
u(t) \,dt \geq
0$, so that the constraint $ ( \int_{\mathbb{R}^d} u(t) \,dt
 )^2 \leq
1$ becomes $ \int_{\mathbb{R}^d} u(t) \,dt \leq1$. In this situation, the
optimization problem addressed in Lemma~\ref{minimisationproblem} can
be solved by variational calculus. However, an explicit form of the
solution is available only if we assume that $u\in\mathcal
{C}^2(B (0,\frac
{R}{2} ))$, meaning that $u$ is twice continuously differentiable
on its support. It is given by the following lemma, which completes the
proof of Proposition~\ref{compactDPP2}.

\begin{lemma}\label{variationalcalculus}
If a function $u$ minimizes $\int_{\mathbb{R}^d} |\nabla u (x) |^2
\,dx$ among
all functions $u$ verifying $\mathcal{M}(\rho,R)$, $u\in\mathcal
{C}^2(B (0,\frac
{R}{2} ))$ and $\int_{\mathbb{R}^d} u(x) \,dx \leq1 $, then $u$
is of the form
\[
u(x) = \biggl(\beta+ \gamma\frac{J_{({d-2})/{2}} (
|x|/\alpha
 )}{|x|^{({d-2})/{2}}} \biggr) \mathbf{1}_{ \lbrace|x|
<{R}/{2}
 \rbrace},
\]
where $\alpha>0,\beta\geq0$ and $\gamma$ are three constants linked
by the conditions $\mathcal{M}(\rho,R)$ and $\int_{\mathbb{R}^d}
u(x) \,dx \leq1 $.
\end{lemma}

\begin{pf*}{Proof of Lemma~\protect\ref{minimisationproblem}}\label
{minimisationproblemproof}
Let $C$ be a kernel which is twice differentiable at $0$ and verifies
the condition $\mathcal{K}_c(\rho,R)$. This implies that $C$ is twice
differentiable everywhere. Moreover, the quantity $\int_{\mathbb
{R}^d} |x|^2 \mathcal{F}
(C)(x)\,dx$ is invariant under radialization of the kernel $C$; see \cite
{ehm2004convolution}, relation (44). Thus, we can consider $C$ as a
radial function. Then, by \cite{ehm2004convolution}, Theorem~3.8,
there exists a countable set $A$ and a sequence of real valued functions
$ \lbrace u _k  \rbrace_{k\in A}$ in $L^2( \mathbb{R}^d)$
such that
\begin{equation}
\label{decompositioncovariancefunction} C(x) = \sum_{k\in A } u_k
\ast u_k (x).
\end{equation}
Further, the convergence of the series is uniform and for each $k\in
A$, the support of $u_k$ lies in $B (0,\frac{R}{2} )$.
Thus,
\begin{equation}
\label{majorationdenablauparfourierdeC} \int_{\mathbb{R}^d} |x|^2 \mathcal{F}(C) (x)
\,dx = \int_{\mathbb
{R}^d} |x|^2 \sum
_{k\in A} \bigl\llvert \mathcal{F} (u_k) (x)\bigr
\rrvert ^2\, dx = \sum_{k\in A}\sum
_{j=1}^{d} \int_{\mathbb
{R}^d}\bigl
\llvert x_j \mathcal{F}(u_k) (x) \bigr\rrvert
^2 \,dx,
\end{equation}
where $x_j$ denotes the $j$th coordinate of the vector $x$. In
addition, we note that $u_k\in L^2(\mathbb{R}^d)$ so $|\cdot|
\mathcal{F}(u_k)(\cdot)
\in L^2(\mathbb{R}^d)$ by~(\ref{majorationdenablauparfourierdeC}). Then,
by \cite{liebanalysis2001}, Theorem~7.9, $\nabla u_k \in L^2(\mathbb{R}^d)$
where $\nabla u_k$ has to be viewed in the distributional sense and
\begin{equation}
\label{derivefouriersobolev} \mathcal{F} ( \partial_j u_k ) (x) = 2i
\pi x_j \mathcal {F}(u_k) (x).
\end{equation}
Thus, from (\ref{majorationdenablauparfourierdeC}) and (\ref
{derivefouriersobolev}) and the Parseval equality,
\[
\int_{\mathbb{R}^d} |x|^2 \mathcal{F}(C) (x)\,dx = \sum
_{k\in A} \int_{\mathbb{R}^d}
\frac{\llvert \nabla u_k (x) \rrvert ^2}{4\pi^2}\,dx.
\]
As every term in the sum above is positive and since this equality
holds for every kernel $C$, the minimum of $\int_{\mathbb{R}^d} |x|^2
\mathcal{F}
(C)(x)\,dx$ is reached if and only if this sum reduces to one term where
$u_k=u$. Then we have $C=u \ast u$ and
\begin{equation}
\label{relationCetu} \int_{\mathbb{R}^d} |x|^2 \mathcal{F}(C) (x)
\,dx = \int_{\mathbb
{R}^d} \frac{\llvert \nabla u (x)
\rrvert ^2}{4\pi^2}\,dx.
\end{equation}
Therefore, minimizing $\int_{\mathbb{R}^d} |x|^2 \mathcal{F}(C)(x)\,
dx$ is equivalent to
minimize $ \int_{\mathbb{R}^d} \llvert \nabla u (x) \rrvert ^2\,dx$.
Hence, it
remains to see what the constraints on the kernel $C$ means for the
function $u$.
Since $C=u \ast u$, where $u$ is one of the function in the
decomposition~(\ref{decompositioncovariancefunction}), $u$ is a
so-called real valued Boas--Kac root of $C$; see \cite
{ehm2004convolution}. Thus, since $C$ is radial, we have by
\cite{ehm2004convolution}, Theorem~3.1, that $u$ is radial and verifies
$u(x)=0$ for $|x|\geq\frac{R}{2}$.
Since $C$ verifies $\mathcal{K}_c(\rho,R)$, we have $C(0)=\rho$ and
$0\leq\mathcal{F}(C)
\leq1$. These constraints are equivalent on $u$ to $\int_{\mathbb{R}^d}
u(x)^2\,dx=\rho$ and $\mathcal{F}(u)^2\leq1$, respectively.
Therefore, $u$~verifies condition $\mathcal{M}(\rho,R)$ and $\mathcal
{F}(u)^2\leq1$.
\end{pf*}

\begin{pf*}{Proof of  Lemma~\protect\ref{existenceminimsobolev}}
According to Lemma~\ref{minimisationproblem}, $C_R$ is a is solution
to Problem~\ref{A} if and only if $C_R=u \ast u$ where $u$ minimizes $\int_{\mathbb{R}
^d} |\nabla u (x) |^2 \,dx$ among all functions $u$ verifying $\mathcal
{M}(\rho
,R)$ and $\mathcal{F}(u)^2\leq1$. We prove the existence of such a
minimum $u$.

Let $\Omega$ denote the open Euclidean ball $B (0,\frac
{R}{2}
)$. Consider the Sobolev space
\[
H^1(\Omega) = \bigl\lbrace f : \Omega\rightarrow\mathbb{R}, f \in
L^2(\Omega), \nabla f \in L^2(\Omega) \bigr\rbrace,
\]
with the norm $\|f\|_{H^1(\Omega)}= ( \|f\|^2 + \| \nabla f \|^2
 )^{\frac{1}{2}}$. For a review on Sobolev spaces, see, for example,
\cite{evans1998PDE} or \cite{liebanalysis2001}. For any $f\in
H^1(\Omega)$, we consider its extension to $\mathbb{R}^d$ by setting $f(x)=0$
if $x\notin\Omega$, so that $f\in L^2(\mathbb{R}^d)$. Let us further denote
$\mathcal{E}$ the set of functions $f\in H^1(\Omega)$ verifying
$\mathcal{M}(\rho
,R)$ and $\mathcal{F}(f)^2\leq1$.

If the minimum $u$ above exists but $u\notin H^1(\Omega)$, then $\int_{\Omega} |\nabla u(x) |^2 \,dx=\infty$, which means that $\mathcal E$ is
empty, otherwise $u$ would not be the solution of our optimization
problem. But $\mathcal{E}$ is not empty (see, e.g., the
functions in Section~\ref{compactfamilysection}), so if $u$ exists,
$u\in H^1(\Omega)$. Let $(w_k)_{k\in\mathbb{N}}$ be a minimizing
sequence in
$\mathcal E$, that is,
\begin{equation}
\label{suiteminimisante} \int_\Omega\bigl\llvert \nabla
w_k(x)\bigr\rrvert ^2 \,dx \mathop{\longrightarrow}_{k\rightarrow+\infty}
\mathop{\inf}_{v
\in\mathcal E} \int_{\Omega
} \bigl\llvert \nabla
v(x)\bigr\rrvert ^2 \,dx,
\end{equation}
where for all $k$, $w_k\in\mathcal E$.
By~(\ref{suiteminimisante}) and since for all $k$, $\int_{\Omega}
|w_k(x)|^2 \,dx = \rho$, the sequence $\lbrace w_k \rbrace$ is bounded in
$H^1(\Omega)$. Then, by the Rellich--Kondrachov compactness theorem
(see \cite{evans1998PDE}), it follows that, up to a subsequence,
$\lbrace w_k \rbrace$ converges in $L^2(\mathbb{R}^d)$ to a certain function
$w\in L^2(\mathbb{R}^d)$ verifying
\begin{equation}
\label{limitesuiteminimisante} \int_\Omega\bigl\llvert \nabla w(x)\bigr\rrvert
^2 \,dx = \mathop{\inf}_{v
\in\mathcal
E} \int_{\Omega}
\bigl\llvert \nabla v(x)\bigr\rrvert ^2 \,dx.
\end{equation}

We now prove that $w\in\mathcal E$, so that $u=w$ is the solution of
our optimization problem. First, $w\in H^1(\Omega)$ as justified
earlier and so $w\in L^2(\mathbb{R}^d)$. Second, as rotations are isometric
functions and since any $w_k$ is radial by hypothesis, we have\vspace*{-2pt} for any
$j\in \mathrm{SO}(d)$
\begin{eqnarray*}
\biggl\lbrace\int_{\mathbb{R}^d} \bigl\llvert w(x) -
w_k(x) \bigr\rrvert ^2 \,dx \rightarrow 0 \biggr\rbrace&
\Longleftrightarrow & \biggl\lbrace\int_{\mathbb{R}^d} \bigl\llvert w
\bigl(j(x)\bigr) - w_k\bigl(j(x)\bigr) \bigr\rrvert ^2
\,dx \rightarrow0 \biggr\rbrace
\\[-2pt]
&\Longleftrightarrow &  \biggl\lbrace\int_{\mathbb{R}^d} \bigl\llvert w
\bigl(j(x)\bigr) - w_k(x) \bigr\rrvert ^2 \,dx
\rightarrow0 \biggr\rbrace.
\end{eqnarray*}
Hence, by uniqueness of the limit, the function $w$ is radial and in
particular, its Fourier transform is real. Further, since $w$ is the
limit in $L^2(\mathbb{R}^d)$ of $w_k$, $w$ verifies the following properties:
\begin{itemize}
\item$w$ is compactly supported in $B ( 0,\frac{R}{2}  )$,
because $w_k\in\mathcal E$ for all $k$.
\item$w\in L^2(\mathbb{R}^d)$ by Rellich--Kondrachov theorem.
\item$\int_{\mathbb{R}^d} |w(x)|^2 \,dx = \int_{\mathbb{R}^d}
|w_k(x)|^2 \,dx =\rho$ since
a sphere in $L^2(\mathbb{R}^d)$ is closed.
\end{itemize}
Therefore, $w$ verifies $\mathcal{M}(\rho,R)$. Third, for every $k$,
$w_k$ being
compactly supported and in $L^2(\mathbb{R}^d)$,
$w_k \in L^1(\mathbb{R}^d)$ so we can consider $\mathcal{F}(w_k)(x)$
for every $x\in\mathbb{R}
^d$ and by the Cauchy--Schwarz inequality\vspace*{-8pt}
\begin{eqnarray*}
&& \bigl\llvert \mathcal{F}(w) (x) - \mathcal{F}(w_k) (x) \bigr
\rrvert \leq a \sqrt{\int_{\mathbb{R}^d} \bigl\llvert
w(t)-w_k(t)\bigr\rrvert ^2 \,dt }\qquad\forall x\in
\mathbb{R}^d,
\end{eqnarray*}
where $a$ is a positive constant. Thereby the convergence of $w_k$ to
$w$ in $L^2(\mathbb{R}^d)$ implies the pointwise convergence of
$\mathcal{F}(w_k)$ to $\mathcal{F}
(w)$. Finally,\vspace*{-2pt} from the relation
\begin{eqnarray*}
&& \mathcal{F}(w_k) (x) \leq1\qquad\forall x\in\mathbb{R}^d,
\forall k\in\mathbb{N},
\end{eqnarray*}
we deduce $\mathcal{F}(w) \leq1$.
\end{pf*}

\begin{pf*}{Proof of Lemma~\protect\ref{variationalcalculus}}
We denote as before $\Omega=B (0,\frac{R}{2} )$.
The optimization problem in Lemma~\ref{variationalcalculus} is a
variational problem with isoperimetric constraints. By \cite
{giaquinta1996variation}, Chapter~2, Theorem~2, every solution must\vspace*{-5pt} solve
\begin{eqnarray}
\Delta u + \lambda_1 u - \frac{ \lambda_2}{2} &=& 0 \qquad
\mbox{on $\Omega $,}
\nonumber
\\[-10pt]
\label{PDEvariationalcalculus}
\\[-10pt]
\nonumber
u &=& 0 \qquad\mbox{on $\partial\Omega$.}
\end{eqnarray}
In equation~(\ref{PDEvariationalcalculus}), $\lambda_1$ and $\lambda
_2$ are the Lagrange multipliers associated to the constraints $ \int
u^2 = \rho$ and $\int u \leq1$, respectively. By the
Karush--Kuhn--Tucker theorem (see~\cite{hiriart1993convex},
Section~VII), \mbox{$\lambda_2 \geq0$}.
Moreover, a solution to the partial differential equation with boundary
condition~(\ref{PDEvariationalcalculus}) is obtained by linear
combination of a homogeneous solution and a particular solution. By
\cite{evans1998PDE}, Section~6.5, Theorem~2, the Laplacian operator
$-\Delta$ has only positive eigenvalues. Hence, the associated
homogeneous equation $ \Delta u + \lambda_1 u=0$ can have a solution
only if $\lambda_1 >0$.

In addition, the function $u$ is radial by hypothesis, so there exists
a function $\tilde{u}$ on $\mathbb{R}$ such that
$u(x)=\tilde{u}( |x| )$ for all $x \in\mathbb{R}^d$.
The partial differential equation~(\ref{PDEvariationalcalculus}) then\vspace*{-2pt} becomes
\begin{eqnarray*}
\tilde{u}''(t) + \frac{d-1}{t}
\tilde{u}'(t) + \lambda_1 \tilde{u}(t) -
\frac{\lambda_2}{2} &=& 0 \qquad\forall t\in \biggl]0,\frac
{R}{2} \biggr[,
\\[-3pt]
\tilde{u} \biggl(\frac{R}{2} \biggr) &=& 0.
\end{eqnarray*}
As $\lambda_1$ is positive, we obtain from \cite{watson1995bessel},
Section~4.31, relations (3) and (4), that a solution to this
equation is of the\vspace*{-3pt} form
\begin{eqnarray}
\label{tildeu} \tilde{u}(t) &=& \biggl(\frac{\lambda_2}{2\lambda_1} + c_1
\frac
{J_{(d-2)/2}(\sqrt{\lambda_1} t)}{t^{(d-2)/2}} + c_2 \frac
{Y_{(d-2)/2}(\sqrt{\lambda_1} t)}{t^{(d-2)/2}} \biggr)
\mathbf{1}_{
 \lbrace
0<t<{R}/2  \rbrace},
\end{eqnarray}
where $Y_{(d-2)/2}$ denotes the Bessel function of the second kind. By
hypothesis, the function $u$ is continuous on $\Omega$ and so at $0$.
Since $Y_{(d-2)/2}$ has a discontinuity at $0$ (see, e.g., \cite
{stegun}) and the remaining terms in (\ref{tildeu}) are continuous, we
must have $c_2=0$. Then, by renaming the constant $c_1$ by $\gamma$ and
letting $\alpha=1/\sqrt{\lambda_1}$, $\beta=\lambda_2/(2\lambda
_1)$, we
obtain that if $u$ is solution to the optimization problem of
Lemma~\ref
{variationalcalculus}, then $u$\vspace*{-5pt} writes
\begin{eqnarray}
u(x) &= & \biggl(\beta+ \gamma\frac{J_{(d-2)/2}(|x|/\alpha
)}{|x|^{(d-2)/2}} \biggr) \mathbf{1}_{ \lbrace x\in\Omega
 \rbrace},
\end{eqnarray}
where $\alpha>0$ and $\beta\geq0$.\vspace*{-4pt}
\end{pf*}

\subsection{Proof of Proposition~\texorpdfstring{\protect\ref{DPPB}}{5.1}}

Let $C$ be given by~(\ref{bessel-type}). According to Proposition~\ref
{DPPexistence}, $\operatorname{DPP}(C)$ exists and has intensity $\rho
$ if $C$
verifies the condition $\mathcal{K}(\rho)$. By \cite{stegun},
equation (9.1.7), we
have $C(0)=\rho$. It is immediate that $C$ is a symmetric continuous
real-valued function. Since Bessel functions are analytic and by the
asymptotic form in \cite{stegun}, (9.2.1), it is clear that $C$
belongs to $L^2(\mathbb{R}^d)$. It remains to obtain $\mathcal F(C)$
and verify
the condition $0 \leq\mathcal{F}( C) \leq1$.

Define\vspace*{-3pt}
\begin{equation}
p_{\sigma}(x) = \frac{J_{({\sigma+d})/{2}}(|x|)}{|x|^{{(\sigma+d)}/{2}}}\qquad\forall x\in\mathbb{R}^d.
\end{equation}
As $p_\sigma$ is radial, by \cite{grafakos2008classical},\vspace*{-3pt} Appendix B.5,
\begin{eqnarray*}
\mathcal{F}(p_{\sigma}) (x) &=& \frac{2\pi}{|x|^{({d-2})/{2}}} \int
_0^{+\infty} r^{{d}/ {2}}
p_{\sigma}(r) J_{({d-2})/{2}}\bigl(2\pi r |x|\bigr)\,dr.
\end{eqnarray*}
By \cite{gradshteyn2007table}, Formula~6.575, we have for\vspace*{-3pt} $\sigma> -2$
\begin{eqnarray*}
&& \mathcal{F}(p_{\sigma}) (x ) = \frac{2\pi}{|2\pi x|^{({d-2})/{2}}} \frac{(1-|2\pi x|^2)_{+}^{{\sigma}/{2}}|2 \pi x|^{({d-2})/{2}}}{2^{{\sigma}/{2}}\Gamma({\sigma}/{2}+1)} =
2^{{d}/{2}-{\sigma}/{2}}\pi^{{d}/{2}}\frac{(1-|2\pi x|^2)^{{\sigma}/{2}}_{+}}{\Gamma(({\sigma+2})/{2})}.
\end{eqnarray*}
Since $C(x)=\rho2^{({\sigma+d})/{2}}\Gamma (\frac{\sigma+d+2}{2} ) p_\sigma (2\frac{x}{\alpha}\sqrt{\frac
{\sigma
+d}{2}}  )$, we obtain~(\ref{Fouriertransformgeneralizedjinc})
by dilatation of the Fourier transform.

We have obviously $\mathcal{F}(C)\geq0$. Since $\sigma\geq0$,
$\mathcal{F}(C)$ attains
its maximum at $0$. Thus, $\mathcal{F}(C)\leq1$ if and only if\vspace*{-4pt}
\[
\mathcal{F}(C) (0)=\frac{\rho(2\pi)^{{d}/{2}}\alpha^d \Gamma
(({\sigma
+d+2})/{2}) }{(\sigma+d)^{{d}/{2}} \Gamma(({\sigma+2})/{2})}\leq1,
\]
which\vspace*{1pt} is equivalent\vadjust{\goodbreak} to $\alpha^d \leq\frac{(\sigma+d)^{{d}/{2}}
\Gamma (({\sigma+2})/{2} )}{\rho(2\pi)^{{d}/{2}}
\Gamma
 (({\sigma+d+2})/{2} ) } $.

Finally, when $\sigma=0$ and $\alpha= \alpha_{\max}$,
$\operatorname{DPP}(C)$ exists
and a straightforward calculation gives $C=C_{B}$.
The convergence result (\ref{Bconvergence}) may be found in \cite
{flyer2006exact} and is a direct application
of~\cite{schoenberg1938metric}, relation (1.8).

\subsection{Proof of Proposition~\texorpdfstring{\protect\ref{DPPlaguerre}}{5.3}}\label
{proofDPPlaguerre}

Define, for all $m\in\mathbb{N}$,
\begin{equation}
f_m (x) = L_m^{d/2} \bigl(|x|^2
\bigr)e^{-|x|^2}\qquad\forall x\in \mathbb{R}^d.
\end{equation}
This function is radial, thus by \cite{grafakos2008classical},
Appendix B.5, we have
\begin{eqnarray*}
\mathcal{F}(f_m) (x) &= & \frac{2\pi}{|x|^{({d-2})/{2}}} \int
_0^{+\infty} r^{{d}/{2}}L_m^{{d}/ {2}}\bigl(r^2\bigr)e^{-r^2}
J_{({d-2})/{2}}(2 \pi r |x|)\,dr.
\end{eqnarray*}
According to \cite{kolbig1996hankel}, we have
\begin{eqnarray*}
\mathcal{F}(f_m) (x) &=& \frac{2\pi}{|x|^{({d-2})/{2}}} \frac
{(-1)^m}{2}
\biggl( \frac{|2\pi x|}{2} \biggr)^{({d-2})/{2}} e^{-{|2\pi x|^2}/{4}}
L_m^{-1-m} \biggl(\frac{|2\pi x|^2}{4} \biggr)
\\
&=&\pi^{{d}/{2}}(-1)^m e^{-|\pi x|^2} \sum
_{k=0}^m\pmatrix{{-1}\cr {m-k}} \frac{(-1)^k |\pi x|^{2k}}{ k!}
\\
&=& \pi^{{d}/{2}}(-1)^{m} e^{-|\pi x|^2} \sum
_{k=0}^m (-1)^{m-k} \frac{(-1)^k |\pi x|^{2k}}{ k!}.
\end{eqnarray*}
Therefore,
\begin{eqnarray*}
\mathcal{F}(f_m) (x) &=& \pi^{{d}/{2}} e^{-|\pi x|^2} \sum
_{k=0}^m \frac{|\pi
x|^{2k}}{ k!}.
\end{eqnarray*}
As $C(x)= \frac{\rho}{{{m-1+{d}/{2}}\choose{m-1}}}
f_{m-1}(\frac
{1}{\sqrt{m}} \frac{x}{\alpha})$, we obtain~(\ref
{Fourierlaguerrefamily}) by dilatation and linearity of the Fourier transform.

Clearly, $\mathcal{F}(C) \geq0$. Thus, we investigate the condition
$\mathcal{F}(C)\leq
1$ for the existence of $\operatorname{DPP}(C)$.
We notice from~(\ref{Fourierlaguerrefamily}) that
\begin{eqnarray}
&&\mathcal{F}(C) (x) = a e^{- b|x|^2} \sum_{k=0}^{m-1}
\frac{b^k|x|^{2k}}{k!},
\end{eqnarray}
where $a$ and $b$ are positive constants. Since $\mathcal{F}(C)$
depends on the
variable $x$ only through its norm, we consider the function $h$ define
for all $r\geq0$ by
$h (r) = \mathcal{F}(C) (  (r,0,\ldots,0 )  )$,
so that for all $x\in\mathbb{R}^d $, $\mathcal{F}(C)(x) = h(|x|)$.
For every $r>0$, $h$ is differentiable at $r$ and a straightforward
calculation leads\vspace*{-3pt} to
\[
h'(r) = a e^{- br^2} \Biggl( -2br \sum
_{k=0}^{m-1} \frac{b^k
r^{2k}}{k!} + \sum
_{k=1}^{m-1}2k \frac{b^k r^{2k-1}}{k!} \Biggr) = - 2 a
e^{- br^2} \frac{b^m r^{2m-1}}{(m-1)!}.
\]
Thus, the function $h$ is decreasing on $(0,+\infty)$. Since $h$ is
continuous on $\mathbb{R}^+$, its maximum is attained at zero, so for every\vspace*{-4pt}
$x\in\mathbb{R}^d$,
\begin{eqnarray*}
&& \mathcal{F}(C) (x) \leq\mathcal{F}(C) (0) = \frac{\rho
(m\pi )^{{d}/{2}}}{{{m-1+{d}/{2}}\choose{m-1}}}
\alpha^d.
\end{eqnarray*}
Hence, $\mathcal{F}(C)\leq1$ if and only if $ \alpha^d \leq\frac
{{{m-1+{d}/{2}}\choose{m-1}}}{\rho(m\pi)^{{d}/{2}}}$.
Moreover $C$ is
radial and since $L_{m-1}^{d/2}(0)={{m-1+{d}/{2}}\choose{m-1}}$, see
\cite{stegun}, relation (22.4.7),\vspace*{1pt} we have $C(0)=\rho$. Therefore, $C$
verifies the condition~$\mathcal{K}(\rho)$ and by Proposition~\ref
{DPPexistence},
$\operatorname{DPP}(C)$ exists and is stationary with intensity $\rho>0$.

It remains to prove the convergence results~(\ref{LGconvergence1})
and~(\ref{LGconvergence2}).
An immediate application of \cite{szego1992orthongonal}, Theorem~8.1.3,
gives the convergence~(\ref{LGconvergence1}), see also~\cite
{alvarez2004laguerre}, Proposition~$1$.\vspace*{-3pt}
Moreover,
\begin{eqnarray}
\label{convergencealphan} &&\lim_{m\rightarrow+\infty} \alpha_{\max} =
\frac{1}{\sqrt{\pi
}\Gamma
 ({d}/{2}+1  )^{{1}/{d}}\rho^{{1}/{d}}}.
\end{eqnarray}
Hence, by~(\ref{LGconvergence1}) and~(\ref{convergencealphan}), we
obtain the\vspace*{-3pt} convergence~(\ref{LGconvergence2}).

\subsection{Proof of Proposition~\texorpdfstring{\protect\ref{DPPcompact}}{5.4}}

By the discussion in Section~\ref{sectionDPPregularitypcffiniterange},
$\operatorname{DPP}(C)$ exists and is an $R$-dependent DPP with
intensity $\rho
$ if $C$ verifies $\mathcal{K}_c(\rho,R)$.
Since $u\in L^2(\mathbb{R}^d)$, the kernel $C$ is continuous by~\cite
{liebanalysis2001}, Theorem~2.20. Moreover, $u(x)=0$ for $|x|>\frac
{R}{2}$, so
by product convolution properties, $C(x)=0$ for $|x|>R$. Hence, $C$
belongs to $L^2(\mathbb{R}^d)$. Since $u$ is radial, so is $C$.
It remains to verify that $0\leq\mathcal{F}(C) \leq1$ and $C(0)=\rho$.

By product convolution properties, we have $C(0)= \int_{\mathbb{R}^d} u(x)^2
\,dx$. From the definition of $u$ in~(\ref{familyu}), we have
\begin{eqnarray*}
&&\frac{\int_{\mathbb{R}^d} u^2(x) \,dx}{\rho\beta(R,\alpha)^2}
\\[-2pt]
&&\quad=\int_{\mathbb{R}^d} \biggl( 1- 2 \biggl(\frac{R}{2}
\biggr)^{{d}/{2}-1} \frac
{J_{({d-2})/{2}}(\llvert {x}/{\alpha}\rrvert )}{J_{({d-2})/{2}}
({R}/({2\alpha})) |x|^{({d-2})/{2}}}\\
&&\qquad {}+ \biggl(\frac{R}{2}\biggr)^{d-2}
\frac{J^2_{({d-2})/{2}}(\llvert {x}/{\alpha}\rrvert )}{J^2_{({d-2})/{2}}({R}/({2\alpha}))|x|^{d-2}} \biggr) \mathbf{1}_{
\lbrace|x|\leq
{R}/{2}  \rbrace} \,dx
\\[-2pt]
&&\quad=\frac{2\pi^{{d}/{2}}}{\Gamma(d/2)} \int_0^{{R}/{2}} \biggl(
r^{d-1} - 2 \biggl( \frac{R}{2} \biggr)^{({d-2})/{2}}
\frac
{J_{({d-2})/{2}}({r}/{\alpha})}{J_{({d-2})/{2}}({R}/({2\alpha}))} r^{{d}/{2}}\\[-2pt]
&&\qquad{}+ \biggl(\frac{R}{2} \biggr)
^{d-2} \frac{J^2_{({d-2})/{2}}({r}/{\alpha})}{J^2_{({d-2})/{2}}({R}/({2\alpha}))} r \biggr) \,dr.
\end{eqnarray*}
By properties of Bessel functions (see~\cite{stegun}), we notice that
for all $b\in\mathbb{R}$, a primitive of $xJ^2_{({d-2})/{2}}(bx)$ is
given by $
\frac{x^2}{2}  (J_{({d-2})/{2}}^2(xb)-J_{{{d}/{2}-2}}(xb)J_{{d}/{2}}(xb)  )$. It\vspace*{1.5pt} follows from~\cite
{grafakos2008classical}, Appendix~B.3, that
\begin{eqnarray*}
\frac{\int_{\mathbb{R}^d} u^2(x)\,dx}{\rho\beta(R,\alpha)^2} &=& \frac{\pi^{d/2}
R^{d}}{d \Gamma({d}/{2})2^{d-1}} - 4 \biggl(\frac{R}{2}
\biggr)^{d-1} \frac{\pi^{{d}/{2}}\alpha}{\Gamma({d}/{2})} \frac
{J_{{d}/{2}}({R}/({2\alpha}))}{J_{{d}/{2}-1}({R}/({2\alpha}))}
\\
&&{}+ \biggl(\frac{R}{2} \biggr)^{d}\frac{\pi^{{d}/{2}}}{\Gamma
({d}/{2})} \biggl( 1
- \frac{J_{{d}/{2}-2}({R}/({2\alpha}))J_{{d}/{2}}({R}/({2\alpha}))}{J_{{d}/{2}-1}^2({R}/({2\alpha}))} \biggr).
\end{eqnarray*}
Thus,\vspace*{1pt} by the definition of $\beta(R,\alpha)$, we obtain that $\int_{\mathbb{R}
^d} u(x)^2\,dx = \rho$.

We now calculate $\mathcal F(C)$. We have $\mathcal{F}(C)= \mathcal
{F}(u)^2$. Since $u$
is radial, $\mathcal{F}(u)$ is real valued and so $\mathcal{F}(C)\geq
0$. In addition, we
have by \cite{grafakos2008classical}, Appendix B.5 and~(\ref{familyu}),
\begin{eqnarray*}
\mathcal{F}(u) (x) &=&\sqrt{\rho} \beta(R,\alpha) \frac{2\pi
}{|x|^{({d-2})/{2}}}\biggl(
\int_0^{{R}/{2}} r^{{d}/{2}}
J_{({d-2})/{2}}\bigl(2\pi r |x|\bigr)\,dr
\\
&&{}- \frac{R^{{d}/{2}-1}}{
2^{{d}/{2}-1} J_{{d}/{2}-1}({R}/({2\alpha}))} \int_0^{{R}/{2}} r
J_{({d-2})/{2}} \biggl(\frac{r}{\alpha} \biggr) J_{({d-2})/{2}}\bigl(2\pi r |x|\bigr)
\,dr \biggr).
\end{eqnarray*}
Since $\alpha>0$, we have by \cite{grafakos2008classical},
Appendix~B.3 and \cite{gradshteyn2007table}, formula~(6.521),
\begin{eqnarray}
&&\!\!\mathcal{F}(u) (x)\nonumber\\
&&\!\!\quad= \sqrt{\rho} \beta(R,\alpha) \frac{2\pi}{|x|^{({d-2})/{2}}} \biggl(
\frac{R^{{d}/{2}} }{\pi2^{{{d}/{2}+1} }}\frac{J_{{d}/{2}}(\pi R |x|)}{|x|}+ \frac{R^{{d}/{2}-1}}{2^{{d}/{2}}J_{{d}/{2}-1}({R}/({2\alpha}))}
\\
\nonumber
&&\!\!\qquad{}\times
\frac{R\alpha J_{({d-2})/{2}}' ({R}/({2\alpha})) J_{({d-2})/{2}}(
\pi
R |x|) - 2\pi\alpha^2 R J_{({d-2})/{2}}( {R}/({2\alpha})) |x|
J_{({d-2})/{2}}'(\pi R |x|) } {1 - 4 \pi^2 |\alpha x|^2 } \biggr)
\end{eqnarray}
from which we deduce the Fourier transform of $u$ in Proposition~\ref
{DPPcompact}. Therefore, if $\alpha$ is such that $\mathcal{F}(u)^2
\leq1$,
then $\mathcal{F}(C)\leq1$ and so $C$ verifies $\mathcal{K}_c(\rho,R)$.

\section*{Acknowledgments}
The authors are grateful to Jean-Fran\c cois Coeurjolly for
illuminating comments and to anonymous referees for numerous
suggestions and comments which helped to improve this paper.

\begin{supplement}
\stitle{Supplement to ``Quantifying repulsiveness of
determinantal point processes''}
\slink[doi]{10.3150/15-BEJ718SUPP} 
\sdatatype{.pdf}
\sfilename{BEJ718\_supp.pdf}
\sdescription{We provide some illustrations of the nonuniqueness of
the most globally repulsive DPP in the sense of Definition~\ref
{globallyrepulsive}, as stated in Proposition~\ref{jinc}. We also show
the p.c.f.s and some realizations associated to different values of the
parameters for the parametric families of DPPs introduced in
Section~\ref{DPPkernelsfamily}.}
\end{supplement}

%





\printhistory
\end{document}